\documentclass[10pt,a4paper]{article}
\pdfoutput=1 

\usepackage{lineno,hyperref}
\modulolinenumbers[5]

\usepackage{amssymb}
\usepackage{amsmath}
\usepackage{algorithm}     
\usepackage{algpseudocode} 
\usepackage{subfig}
\usepackage{csquotes}
\usepackage{tikz}
\usetikzlibrary{arrows}
\usepackage{mathtools}
\usepackage{bbm}
\usepackage{amsthm}
\usepackage{pgfplots}
\usepackage{float}
\usepgfplotslibrary{external}

\newcommand{\email}[1]{\hspace*{\stretch{1}}\emph{\texttt{#1}}}

\makeatletter
\def\blfootnote{\xdef\@thefnmark{$\star$}\@footnotetext}
\makeatother
\newenvironment{Authors}%
  {\begin{center}\begin{bfseries}}%
  {\end{bfseries}\end{center}}
\newenvironment{Addresses}%
  {\begin{flushleft}\begin{itshape}}%
  {\end{itshape}\end{flushleft}}

 \usepackage{fancyhdr}  
  
\fancypagestyle{plain}{
	\fancyhead{}
	\fancyhead[C]{\hfill submitted to Journal of Scientific Computing, January 2021}
}
  
\newtheorem{definition}{definition}[section]
\newtheorem{theorem}{Theorem}[section]
\newtheorem{theorem1}{Theorem}[section]
\newtheorem{theorem2}{Theorem}[section]
\newtheorem{proposition}[theorem]{Proposition}
\newtheorem{lemma}[theorem1]{Lemma}
\newtheorem{remark}[theorem2]{Remark}

  \newcommand{\vertiii}[1]{{\left\vert\kern-0.25ex\left\vert\kern-0.25ex\left\vert #1 
    \right\vert\kern-0.25ex\right\vert\kern-0.25ex\right\vert}}

\begin{document}

\thispagestyle{plain}

\title{Registration-based model reduction in complex  two-dimensional geometries.}
 \date{}
 
 \maketitle
\vspace{-50pt} 
 
\begin{Authors}
Tommaso Taddei$^{1}$,
Lei Zhang$^{1}$
\end{Authors}

\begin{Addresses}
$^1$
IMB, UMR 5251, Univ. Bordeaux;  33400, Talence, France.
Inria Bordeaux Sud-Ouest, Team MEMPHIS;  33400, Talence, France, \email{tommaso.taddei@inria.fr,lei.a.zhang@inria.fr} \\
\end{Addresses}

\begin{abstract}
We present  a general — i.e., independent of the underlying equation — registration procedure  for parameterized model order reduction.
Given the spatial domain $\Omega \subset \mathbb{R}^2$ and the manifold
$\mathcal{M}= \{ u_{\mu} : \mu \in \mathcal{P}  \}$ associated with the parameter domain $\mathcal{P} \subset \mathbb{R}^P$ and the parametric field
$\mu \mapsto u_{\mu} \in   L^2(\Omega)$,
our approach takes as input a  set of snapshots $\{ u^k \}_{k=1}^{n_{\rm train}} \subset \mathcal{M}$ and returns a parameter-dependent bijective mapping 
${\Phi}: \Omega \times \mathcal{P} \to \mathbb{R}^2$:
the mapping is designed to make the mapped manifold $\{ u_{\mu} \circ {\Phi}_{\mu}: \, \mu \in \mathcal{P} \}$ more amenable  for linear compression methods. 
In this work, we extend and further analyze the registration approach proposed in 
[Taddei, SISC, 2020].
The contributions of the present work are twofold. 
First, we extend the approach to deal with 
annular domains by introducing a suitable transformation of the coordinate system.
Second, we discuss the extension to 
general two-dimensional geometries: towards this end, we introduce a spectral element  approximation, which relies on a partition  $\{  \Omega_{q} \}_{q=1}
^{N_{\rm dd}}$  
of the domain $\Omega$
such that 
$\Omega_1,\ldots,\Omega_{N_{\rm dd}}$ are 
isomorphic to the unit square.
We further show that our spectral element approximation can cope with  parameterized geometries.
We present rigorous mathematical analysis to justify our proposal; furthermore, 
we present numerical results for 
a heat-transfer problem in an annular domain, 
a    potential flow past a rotating symmetric  airfoil, {and an inviscid transonic compressible flow  past a non-symmetric airfoil}, to demonstrate the effectiveness  of our method.
\end{abstract}


\section{Introduction}
\label{sec:introduction}
\subsection{Registration-based model order reduction}

The inadequacy of linear approximation methods to deal with parametric fields with sharp gradients hinders the application of parameterized model order reduction (pMOR, \cite{hesthaven2016certified,quarteroni2015reduced,rozza2007reduced}) techniques to a broad class of problems, including high-Reynolds flows, contact problems, etc. To address this issue, several authors have proposed to resort to nonlinear approximation methods. The goal of this paper is to develop a \emph{general} (i.e., independent of the underlying model) registration-based data compression technique for problems with slowly-decaying Kolmogorov $N$-widths, \cite{pinkus2012n}; more in detail, we wish to extend the approach introduced in 
\cite{taddei2020registration,taddei2021space} to more general two-dimensional geometries.

We denote by $\mu$ the vector of model parameters in the parameter region $\mathcal{P} \subset \mathbb{R}^P$; we denote by $\Omega \subset \mathbb{R}^2$ the domain of interest, and we denote by $u_{\mu}$ the solution to the   partial differential equation (PDE) of interest for the parameter $\mu \in \mathcal{P}$. 
We define  the Hilbert  space $\mathcal{X}$ equipped with the inner product $(\cdot,\cdot)$ and the induced norm $\| \cdot \| := \sqrt{(\cdot,\cdot)}$; then, we introduce 
the solution manifold $\mathcal{M} = \{ u_{\mu}: \mu \in \mathcal{P}  \} \subset \mathcal{X}$ . 
We here assess our methodology through the vehicle of 
a heat-transfer problem in an annular domain,
a    potential flow past an airfoil, and
a transonic inviscid flow past an airfoil.
We shall consider 
$\mathcal{X} =  H^1(\Omega)$ with
$(w,v) := \int_{\Omega} 
\nabla w \cdot  \nabla v
\, +  \,  w \cdot v \, d{x}$  in the first two examples and 
$\mathcal{X} =  L^2(\Omega)$ with
$(w,v) := \int_{\Omega}  \,  w \cdot v \, d{x}$
in the third example.

The key feature of registration-based (or Lagrangian) 
nonlinear compression methods (e.g., \cite{iollo2014advection,ohlberger2013nonlinear,taddei2015reduced}) is the introduction of  a parametric mapping ${\Phi}: \Omega \times \mathcal{P} \to \Omega$ such that (i) ${\Phi}_{\mu}$ is a bijection from $\Omega$ in itself for all $\mu \in \mathcal{P}$, and 
(ii) the mapped manifold
$\widetilde{\mathcal{M}} = \{ \widetilde{{u}}_{\mu} : = {u}_{\mu}\circ {\Phi}_{\mu} \, :  \, \mu \in \mathcal{P}\}$ is more amenable  for linear  compression methods.  Note that 
mappings have been  broadly used in the pMOR literature to deal with parameterized geometries (see 
\cite{lassila2014model,manzoni2012model} and the references therein); however, as discussed in 
\cite{taddei2020registration}, the use of mappings is here motivated by approximation considerations rather than by the need to restate the PDE on a reference (parameter-independent) configuration.

In  \cite{taddei2021space}, we   developed a general procedure for the construction of affine maps ${\Phi}$ of the form
\begin{equation}
\label{eq:affine_maps}
{\Phi}  = \texttt{id} + W_M \mathbf{a},
\end{equation}
 where $\texttt{id}: \Omega \to \Omega$ is the identity map, and 
 $W_M: \mathbb{R}^M \to   {\rm Lip}(\Omega; \mathbb{R}^2)$ is a suitable linear operator.
Given snapshots of the manifold $\{ {u}^k = u_{\mu^k} \}_{k=1}^{n_{\rm train}} \subset \mathcal{M}$, the approach returns (i) a $N$-dimensional linear operator $Z_N: \mathbb{R}^N \to \mathcal{X}$,
(ii) the  $M$-dimensional linear operator  $W_M$ in \eqref{eq:affine_maps}, and 
 (iii) coefficients
 $\{ \boldsymbol{\alpha}^k \}_{k=1}^{n_{\rm train}} \subset \mathbb{R}^N$ and 
  $\{ \mathbf{a}^k \}_{k=1}^{n_{\rm train}} \subset \mathbb{R}^M$ such that
  \begin{equation}
  \label{eq:output_rePOD_affine}
  u^k \circ {\Phi}^k \approx \widehat{u}^k,
  \quad
  {\Phi}^k = \texttt{id} + W_M \mathbf{a}^k,
  \;\;
  \widehat{u}^k = Z_N  \boldsymbol{\alpha}^k ,
  \quad
  k=1,\ldots,n_{\rm train},
  \end{equation}
 The approach relies on repeated solutions to a non-convex optimization problem to build the mappings $\{ {\Phi}^k \}_k$, and on  proper orthogonal decomposition (POD, \cite{berkooz1993proper,volkwein2011model}) to generate the low-dimensional approximation operators $Z_N, W_M$. 
The approach was successfully applied to the space-time approximation of one-dimensional hyperbolic conservation laws. One major limitation of the approach is the inability to deal with domains that are not isomorphic to the unit square: this limitation precludes the application of the approach in  \cite{taddei2021space}  to two-dimensional steady and unsteady PDEs in general domains.

\subsection{Objective and layout of the paper}

The contributions of the present work are twofold.
\begin{enumerate}
\item
We extend the  approach in \cite{taddei2021space} to annular domains $\Omega = \mathcal{B}_R({0}) \setminus \mathcal{B}_r({0})$ with $0<r<R$. Towards this end, we introduce a polar transformation
${\Psi}: \widehat{\Omega}_{\rm pol} \to \Omega$, $\widehat{\Omega}_{\rm pol} = (0,1) \times (-1/2,1/2)$, such that ${\Psi}(\rho,\theta) = (r + (R-r) \rho) [\cos (2 \pi \theta), \sin (2 \pi \theta)]$ and we consider mappings of the form:
\begin{equation}
\label{eq:spectral_map_annulus}
{\Phi} 
=
{\mathcal{N}}(\cdot,  W_M \mathbf{a})
=
{\Psi} \circ {\Phi}_{\rm pol}  \circ {\Lambda},
\end{equation}
where ${\Phi}_{\rm pol}  = \texttt{id} + W_M \mathbf{a}$, 
 $W_M: \mathbb{R}^M \to  {\rm Lip}(\widehat{\Omega}_{\rm pol} ; \mathbb{R}^2   )$ is a suitable linear operator and ${\Lambda} = {\Psi}^{-1}$.
\item
We extend  the approach to   arbitrary two-dimensional domains through partitioning.
Given the partition $\{ \Omega_q \}_{q=1}^{N_{\rm dd}}$ of $\Omega$, 
we denote by ${\Psi}_q: \widehat{\Omega}  \to  \Omega_q$ the bijective mapping between $\widehat{\Omega}  = (0,1)^2$ and the $q$-th element of the partition,  and we denote by
${\Lambda}_q:     \Omega_q \to \widehat{\Omega}$ the inverse of ${\Psi}_q$. Then, we consider piecewise-smooth mappings of the form
\begin{equation}
\label{eq:piecewise_affine_mappings}
{\Phi} =
{\mathcal{N}}(\cdot,  W_M \mathbf{a})
=
\sum_{q=1}^{N_{\rm dd}} \; \left(
{\Psi}_q \circ  {{\Phi}}_q \circ {\Lambda}_q
\right)
\mathbbm{1}_{\Omega_q},
\end{equation} 
where 
$ {{\Phi}}_q= \texttt{id} + W_M^q \mathbf{a}: \widehat{\Omega}\to \widehat{\Omega}$,
 $\mathbbm{1}_{\Omega_q}$ denotes the indicator function associated with the $q$-th element of the partition,
 and  $
W_M = [W_M^1, \ldots,W_M^{N_{\rm dd}}]: \mathbb{R}^M \to [ {\rm Lip}(\widehat{\Omega}, \mathbb{R}^2)]^{N_{\rm dd}}$ is a  suitable linear operator.   We also discuss how to adapt \eqref{eq:piecewise_affine_mappings} to deal with parameterized domains.
\end{enumerate}

The outline of the paper is as follows.
In section \ref{sec:spectral_maps}, we introduce the approximation spaces employed for registration: we review the case of rectangular domains, and then we discuss in detail the definition of mappings of the forms \eqref{eq:spectral_map_annulus} and \eqref{eq:piecewise_affine_mappings}.
Then, in section \ref{sec:registration}, we present the registration algorithm, while in section \ref{sec:numerics} we present results for the two model problems considered in this work, to demonstrate the effectiveness of our proposals.
Section \ref{sec:conclusions}  completes  the paper.

\subsection{Relationship with previous works}
Registration-based techniques are tightly linked to a number of techniques in  related fields.
First, registration is central in image processing: in this field, registration refers to the process of transforming different sets of data into one coordinate system, \cite{zitova2003image}.
In computational mechanics,
 Persson and Zahr have proposed in \cite{zahr2018optimization}  an $r$-adaptive optimization-based 
high-order discretization method 
 to deal with shocks/sharp gradients of the solutions to advection-dominated problems.
In uncertainty quantification, several authors (see, e.g., \cite{marzouk2016sampling}) have proposed measure transport approaches to sampling: transport maps are used to ``push forward'' samples from a reference configuration and ultimately facilitate sampling from non-Gaussian distributions.
As recently observed in \cite{mojgani2020physics}, the notion of registration is also at the core of diffeomorphic dimensionality reduction (\cite{walder2008diffeomorphic}) in the field of machine learning.

The  partitioned approach  proposed in this paper
is similar in scope to the ``reduced basis triangulation" in \cite{rozza2007reduced}, and  shares important features with  isoparametric spectral element discretizations of PDEs,
\cite{korczak1986isoparametric}.  As in \cite{lovgren2006reduced},  we build the local maps $\{ {\Psi}_q \}_q$ in \eqref{eq:piecewise_affine_mappings}   using Gordon-Hall transformations  \cite{gordon1973construction}.
Furthermore, the adaptive construction of the parameterized
maps $\{ {\Phi}_q \}_q$ in the unit square relies on several   building blocks of the original proposals in  \cite{taddei2020registration,taddei2021space}. 
Finally, enforcement of discrete bijectivity 
 (cf. Definition \ref{def:discrete_bijectivity})
exploits the same mesh distortion indicator employed in \cite{zahr2018optimization}.

In recent years,  there has been a growing interest in the development of nonlinear reduction techniques for dimensionality reduction in pMOR.
Several authors have proposed approximations of the form 
$\widehat{u}_{\mu} = Z_{N,\mu} ( \widehat{\boldsymbol{\alpha}}_{\mu} )$ where
$Z_{N}: \mathbb{R}^N \times \mathcal{P} \to \mathcal{X}$ is a nonlinear  and/or parameter-dependent operator, which is built based on 
Grassmannian learning \cite{amsallem2008interpolation,zimmermann2018geometric},
convolutional auto-encoders \cite{fresca2021comprehensive,kashima2016nonlinear,kim2020efficient,lee2020model},
transported/transformed snapshot methods
\cite{cagniart2019model,nair2018transported,reiss2018shifted,welper2017interpolation},
 displacement interpolation \cite{rim2018model}.
In this paper, we do not discuss in detail these methods and their relation with registration-based techniques. 

The ultimate goal of pMOR is to exploit the results of   data compression techniques --- here, 
the operators $Z_N, W_M$ and the solution and mapping coefficients   $\{ \boldsymbol{\alpha}^k \}_k \subset \mathbb{R}^N$ and 
  $\{ \mathbf{a}^k \}_k\subset \mathbb{R}^M$ --- 
  to estimate the solution field ${u}_{\mu}$ for new values of the parameter $\mu$ in $\mathcal{P}$.
In this work, we pursue a fully non-intrusive
(\cite{chakir2009methode,gallinari2018reduced,guo2018reduced})
 technique 
based on    radial basis function (RBF, \cite{wendland2004scattered}) approximation, 
  to estimate solution and mapping coefficients;
  in \cite{taddei2020registration,taddei2021space}, we resorted to RBF approximation to estimate the mapping coefficients,  and to Galerkin/Petrov-Galerkin projection to estimate the solution coefficients. Since the emphasis of this work is on the treatment of complex geometries, we here choose to not discuss in detail this aspect.

In our previous works \cite{taddei2020registration,taddei2021space} and in the previous section (cf. \eqref{eq:output_rePOD_affine}), we presented  our registration technique as a complete data compression algorithm, which takes as input a set of snapshots, and returns the operators  $Z_N, W_M$ and the coefficients 
 $\{ \boldsymbol{\alpha}^k \}_{k=1}^{n_{\rm train}} \subset \mathbb{R}^N$,
  $\{ \mathbf{a}^k \}_{k=1}^{n_{\rm train}} \subset \mathbb{R}^M$ such that
$u^k \circ {\Phi}^k \approx \widehat{u}^k$,
${\Phi}^k =  {\mathcal{N}}(\cdot,  W_M \mathbf{a}^k)$, 
$\widehat{u}^k = Z_N  \boldsymbol{\alpha}^k$ ,
$k=1,\ldots,n_{\rm train}$.
In this paper, we shall interpret registration based on Algorithm \ref{alg:registration} and on the subsequent generalization (cf. \eqref{eq:parametric_mapping_Phi}) as a preliminary preconditioning step that ``simplifies'' --- in the sense of data compression --- the task of model reduction. 
From an implementation standpoint, registration is performed \emph{before} applying (either intrusive or non-intrusive) model reduction:  it can thus be easily integrated with existing pMOR routines for parameterized geometries.
Furthermore,  we do not have to consider the same training parameters $\{  \mu^k \}_k$ to construct the mapping and to construct the mapped solution.

\section{Spectral maps for registration}
\label{sec:spectral_maps}
In this section, we introduce the approximation spaces employed in section \ref{sec:registration} for registration; furthermore, we present rigorous mathematical analysis that provides sufficient and computationally-feasible conditions for the bijectivity of the mapping ${\Phi}$.
In section \ref{sec:registration_approx}, we illustrate how to approximate parametric fields using registration-based model reduction: to fix notation, we assume that the underlying high-fidelity (hf) discretization is based on the finite element (FE) method.
In section \ref{sec:affine_maps} 
we review the special  case of rectangular domains, while in 
section \ref{sec:affine_maps_annulus} we consider annular domains; then, in  section \ref{sec:piecewise_maps}  we address the general case of two-dimensional domains.  In sections \ref{sec:affine_maps} and \ref{sec:affine_maps_annulus}, we consider  spectral affine transformations of the form \eqref{eq:affine_maps} and \eqref{eq:spectral_map_annulus}; in section \ref{sec:piecewise_maps}, we consider spectral element approximations of the form \eqref{eq:piecewise_affine_mappings}.
Finally, in section \ref{sec:bijectivity}, we discuss the practical enforcement of the bijectivity condition  at the continuous and discrete (cf. Definition \ref{def:discrete_bijectivity}) level.

Given the mapping ${\Phi}: \Omega \to \mathbb{R}^2$, we denote by  $g_{\Phi}$ the Jacobian determinant,  $g_{\Phi}  = 
\,{\rm det} \left( {\nabla} {\Phi} \right)$;
given the tensorized two-dimensional domain $U = U_1 \times U_2$, $\mathbb{Q}_J(U)$ refers to the space of tensorized polynomials of degree lower or equal  to $J$,   $\mathbb{Q}_J(U) = {\rm span}\{ p_1(x_1) p_2(x_2) : p_i \in \mathbb{P}_J(U_i),\;i=1,2 \}$, for some $J\geq1$.

\subsection{Approximation of parametric fields using registration}
\label{sec:registration_approx}

\subsubsection{Finite element discretization}
Similarly to \cite{taddei2020discretize}, we consider a FE isoparametric discretization of degree $\texttt{p}$.
Given the domain $\Omega \subset \mathbb{R}^2$,
we define the triangulation $\{ \texttt{D}_k  \}_{k=1}^{N_{\rm e}}$,  where $\texttt{D}_k \subset \Omega$ denotes  the $k$-th element of the mesh.
We define the reference element 
$\widehat{\texttt{D}} = \{ {X}\in (0,1)^2:  X_1+X_2 <  1  \}$  and the bijection  ${\Psi}_k^{\rm hf}$ from  $\widehat{\texttt{D}}$ to $\texttt{D}_k$ for $k=1,\ldots,N_{\rm e}$.
 We define the Lagrangian basis $\{ \ell_i \}_{i=1}^{n_{\rm lp}}$  of the polynomial space $\mathbb{P}_{\texttt{p}}(\widehat{\texttt{D}})$ associated with the nodes 
 $\{  {X}_i  \}_{i=1}^{n_{\rm lp}}$;
 then, we introduce the  mappings $\{ {\Psi}_k^{\rm hf} \}_k$ such that
  \begin{equation}
\label{eq:psi_mapping}
 {\Psi}_k^{\rm hf}( {{X}})
 =
 \sum_{i=1}^{n_{\rm lp}} \;
 {x}_{i,k}^{\rm hf}  \; \ell_i({X}),
\end{equation}
where $\{ {x}_{i,k}^{\rm hf}   := {\Psi}_k^{\rm hf}( {{X}}_i): \, i=1,\ldots,n_{\rm lp}, k=1,\ldots,N_{\rm e}  \}$ are the nodes of the mesh. We define the basis functions
$\ell_{i,k} :=  \ell_i \circ  ( {\Psi}_k^{\rm hf} )^{-1}: \texttt{D}_k \to \mathbb{R}$.  Note that ${\Psi}_k^{\rm hf}$ is completely characterized by  the nodes
in the $k$-th element 
 ${\texttt{X}}_{k}^{\rm hf} :=   \{  {x}_{i,k}^{\rm hf}  \}_{i=1}^{n_{\rm lp}}$, $k=1,\ldots,N_{\rm e}$.
We further introduce the
nodes of the mesh $\{ {x}_j^{\rm hf}  \}_{j=1}^{N_{\rm hf}}$ taken without repetitions  and the
 connectivity matrix $\texttt{T} \in \mathbb{N}^{n_{\rm lp}, N_{\rm e}}$,  such that
${x}_{i,k}^{\rm hf}   = {x}_{  \texttt{T}_{i,k} }^{\rm hf} $. 

\begin{definition}
\label{def:mesh}
\textbf{High-order FE mesh.}
If we fix the reference element  $\widehat{\texttt{D}}$ and the reference nodes 
$\{  {X}_i  \}_{i=1}^{n_{\rm lp}} \subset \overline{\widehat{\texttt{D}}}$, a FE mesh $\mathcal{T}_{\rm hf}$ of degree \texttt{p} of $\Omega$ is uniquely identified by  
the nodes $\{   {x}_j^{\rm hf}\}_{j=1}^{N_{\rm hf}}$ and the  
 connectivity matrix $\texttt{T} \in \mathbb{N}^{n_{\rm lp}, N_{\rm e}}$,
 $\mathcal{T}_{\rm hf} = \left(   \{   {x}_j^{\rm hf}\}_{j=1}^{N_{\rm hf}}, \texttt{T}   \right)$.
 The FE space of order \texttt{p} associated with the mesh $\mathcal{T}_{\rm hf}$ is then given by
\begin{equation}
\label{eq:FE_space}
\mathfrak{X}_{\mathcal{T}_{\rm hf}} : = \left\{
v\in C(\Omega) : \; \;
v    \circ {\Psi}_k^{\rm hf}  \in \mathbb{P}_{\texttt{p}}(\widehat{\texttt{D}}),
\;
k=1,\ldots,N_{\rm e} 
\right\}.
\end{equation}
If $u \in \mathfrak{X}_{\mathcal{T}_{\rm hf}}$, we denote by $\mathbf{u} \in \mathbb{R}^{N_{\rm hf}}$ the vector such that
$( \mathbf{u} )_j = u(  {x}_j^{\rm hf} )$ for $j=1,\ldots,N_{\rm hf}$;
note that
  \begin{equation}
  \label{eq:vector2field}
  u \Big|_{  \texttt{D}_{k}  }  \, = \,  
 \sum_{i=1}^{n_{\rm lp}} \; 
 \left( \mathbf{u}  \right)_{\texttt{T}_{i,k}} \; \,  \ell_{i,k},
 \qquad
 k=1,\ldots,N_{\rm e}.
 \end{equation}
 We further denote by $\mathbf{X} \in \mathbb{R}^{N_{\rm hf}, N_{\rm hf}}$ the symmetric positive definite matrix associated with the inner product  $(\cdot,\cdot)$:
 \begin{equation}
 \label{eq:Xnorm}
 (u, v) = 
 \mathbf{v}^T \, \mathbf{X} \,  \mathbf{u},
 \quad
 \forall \, u,v \in \mathfrak{X}_{\mathcal{T}_{\rm hf}}.
 \end{equation}
\end{definition}

Given the mesh $\mathcal{T}_{\rm hf}$ over $\Omega$ and the 
bijection ${\Phi}: \Omega \to {\Phi} ( \Omega )$, we  define the mapped mesh
 ${\Phi} (\mathcal{T}_{\rm hf}) $ that shares with $\mathcal{T}_{\rm hf}$ the same connectivity matrix $\texttt{T}$ and has nodes $\{ {\Phi} ({x}_j^{\rm hf})    \}_{j=1}^{N_{\rm hf}}$, 
 ${\Phi} (\mathcal{T}_{\rm hf}) = \left( \{ {\Phi} ({x}_j^{\rm hf})    \}_{j=1}^{N_{\rm hf}}, \texttt{T} \right)$. 
 We denote by  $ {\Psi}_{k,\Phi}^{\rm hf}$ the elemental 
 mapping associated with the $k$-th element 
$\texttt{D}_{k,\Phi}$ 
 of  ${\Phi} (\mathcal{T}_{\rm hf})$; $ {\Psi}_{k,\Phi}^{\rm hf}$
  is given by
  \begin{equation}
\label{eq:psi_mapping_mapped}
{\Psi}_{k,\Phi}^{\rm hf}( {{X}})
 =
 \sum_{i=1}^{n_{\rm lp}} \;
 {\Phi}({x}_{i,k}^{\rm hf}  ) \; \ell_i({X}).
\end{equation}
Next Definition    is key for the discussion.

\begin{definition}
\label{def:discrete_bijectivity}
\textbf{Discrete bijectivity.}
Given the mesh $\mathcal{T}_{\rm hf}$, we say that the transformation
${\Phi}: \Omega \to {\Phi} ( \Omega )$ is  bijective for $\mathcal{T}_{\rm hf}$, if the FE mappings $\{ {\Psi}_{k,\Phi}^{\rm hf}  \}_k$ are invertible.
\end{definition}

It is possible to verify that  a 
bijective mapping ${\Phi}$ in $\Omega$ might not satisfy 
discrete bijectivity for a given  mesh $\mathcal{T}_{\rm hf}$; in particular, for highly anisotropic meshes, discrete bijectivity must be explicitly enforced. If we are ultimately interested in performing (Petrov-)Galerkin projection --- as in projection-based pMOR --- we shall ensure that our geometric parameterization ${\Phi}$ 
leads to consistent mapped FE meshes.

\subsubsection{Approximation of parametric fields}
\label{sec:approx_form_gen}

Given the manifold $\mathcal{M} = \{ u_{\mu}: \mu \in \mathcal{P}  \}$, we wish to construct a low-rank approximation of the elements of $\mathcal{M}$ that can be rapidly queried for any $\mu \in \mathcal{P}$. Given the mesh $\mathcal{T}_{\rm hf}$, the parametric mapping ${\Phi}: \Omega \times \mathcal{P} \to \mathbb{R}^2$, the  matrix $\mathbf{Z}_N \in \mathbb{R}^{N_{\rm hf}, N}$ and the function $\widehat{\boldsymbol{\alpha}}: \mathcal{P} \to \mathbb{R}^N$,   registration-based methods  aim to devise  approximations of the form
\begin{equation}
\label{eq:registration_fundamental_approx}
\mu \mapsto \left( \mathcal{T}_{\rm hf,\mu} = {\Phi}_{\mu}(\mathcal{T}_{\rm hf}), \; \widehat{\mathbf{u}}_{\mu} = \mathbf{Z}_N \widehat{\boldsymbol{\alpha}}_{\mu} \right).
\end{equation}
The pair $\left(  \mathcal{T}_{\rm hf,\mu}, \; \widehat{\mathbf{u}}_{\mu}   \right)$ identifies a 
unique FE field in the space $\mathfrak{X}_{\mathcal{T}_{\rm hf,\mu}}$, see \eqref{eq:FE_space}-\eqref{eq:vector2field}. Note that 
$\mu \mapsto \left( \mathcal{T}_{\rm hf}, \; \widehat{\mathbf{u}}_{\mu} \right)$ can be interpreted as an approximation of the mapped field $u_{\mu} \circ {\Phi}_{\mu}$: the linear operator $Z_N: \mathbb{R}^N \to  \mathfrak{X}_{\mathcal{T}_{\rm hf}}$ associated with $\mathbf{Z}_N$ can thus be interpreted as a reduced-order basis for the elements of the mapped manifold
$\widetilde{\mathcal{M}}$.

\subsection{Spectral  maps in rectangular domains}
\label{sec:affine_maps}

Next result is key for the discussion.

\begin{proposition}
\label{th:application_unit_square}
(\cite[Proposition 2.3]{taddei2020registration})
Let $\Omega$ be a rectangular domain;  consider the mapping 
${\Phi} =  \texttt{id}+ {\varphi}$, where
\begin{equation}
\label{eq:boundary_conditions_mapping}
{\varphi}   \cdot {n} \, \equiv  \, 0 
\quad
{\rm on} \; \partial \Omega,
\end{equation}
where $n$ is the outward normal to $\Omega$.
Then, ${\Phi}$ is bijective from $\Omega$ into itself if 
\begin{equation}
\label{eq:condition_equivalent}
\min_{{x} \in \overline{\Omega}} \,  
g_{\Phi}({x})  \, > \, 0.
\end{equation}
\end{proposition}

Condition \eqref{eq:boundary_conditions_mapping} is easy  to enforce numerically: given the integer $J>1$, we here consider $M$-dimensional spaces   of tensorized polynomials,
\begin{equation}
\label{eq:tensorized_polynomial_space}
\mathcal{W}_M = {\rm span} \{ {\varphi}_m  \}_{m=1}^M \subset \mathcal{W}_{\rm hf} : =
\left\{
{\varphi} \in [\mathbb{Q}_J(\Omega)]^2 \, : \,
{\varphi}   \cdot {n} \big|_{\partial \Omega}\, \equiv  \, 0 
\right\},
\end{equation}
and we define the linear operator $W_M : \mathbb{R}^M \to \mathcal{W}_M$,
$W_M \mathbf{a} = \sum_{m=1}^M ( \mathbf{a} )_m {\varphi}_m$. 
Enforcement of \eqref{eq:condition_equivalent} is more involved: we discuss this issue in   section \ref{sec:bijectivity}. 
We remark that, since  the determinant is  a continuous  function and $g_{\Phi}\equiv 1$ if ${\Phi}=\texttt{id}$, 
for any non-empty linear space $\mathcal{W}_M \subset \mathcal{W}_{\rm hf}$,
there exists a ball of finite radius $r>0$ centered in the origin, $B = \mathcal{B}_r({0})$, $r=r(W_M)$,  such that
${\Phi} =  \texttt{id}+ W_M \mathbf{a}$ is a bijection in $\Omega$ for all $\mathbf{a} \in B$.

Mappings satisfying \eqref{eq:boundary_conditions_mapping} map each edge of $\Omega$ in itself and each corner in itself: since tangential displacement is not necessarily zero at the boundary,  we can  rely on mappings of the form \eqref{eq:output_rePOD_affine} to enforce non-trivial deformations on $\partial \Omega$. 
As proved in the next Lemma,  linear combinations of mappings with finite boundary deformations are not bijections  for non-rectangular domains: this represents a fundamental limitation of affine maps and motivates the discussion of the next two sections.

\begin{lemma}
\label{th:affine_bad}
Consider  $\Omega = \{ {x}  \in \mathbb{R}^2: \| {x} \|_2 < 1 \}$ and consider the mappings ${\Phi}_i = \texttt{id} + {\varphi}_i$ for $i=1,2$. Assume that there exists ${x} \in \partial \Omega$ such that 
${\Phi}_1({x}) \neq 
{\Phi}_2({x})$;  then, ${\Phi}_t =(1-t) {\Phi}_1 + t {\Phi}_2 $ is not a bijection in $\Omega$ for any $t \in (0,1)$. 
\end{lemma} 
 \begin{proof}
If  ${\Phi}_t$ is a bijection in $\Omega$, we must have  ${\Phi}_t(\partial \Omega)  =  \partial \Omega$: therefore,
it suffices to show that ${\Phi}_t({x})$ does not belong to the unit circle. To shorten notation, we omit dependence on ${x}$. We first observe that
 $\| {\Phi}_i \|_2= 1$  for $i=1,2$ and thus we have
 $\|  {\varphi}_i \|_2^2 = - 2\texttt{id} \cdot {\varphi}_i$. Then, we observe that
${\Phi}_t =(1-t) {\Phi}_1 + t {\Phi}_2 
=
{\Phi}_1  + t({\varphi}_2-{\varphi}_1)
$ satisfies
$$
\begin{array}{rl}
\| {\Phi}_t  \|_2^2 = &
\| {\Phi}_1  \|_2^2 +
t^2   \| {\varphi}_2-{\varphi}_1  \|_2^2
+ 2 t
\left( \texttt{id} +  {\varphi}_1 \right) \cdot 
\left( {\varphi}_2 - {\varphi}_1\right) 
\\[3mm]
= &
1 +
t^2   \| {\varphi}_2-{\varphi}_1  \|_2^2
+ 2 t
\left(
\texttt{id}  \cdot  {\varphi}_2 -
\texttt{id}  \cdot  {\varphi}_1 +
{\varphi}_2 \cdot  {\varphi}_1
-
\|  {\varphi}_1 \|_2^2
\right) 
\\[3mm]
= &
1 +
t^2   \| {\varphi}_2-{\varphi}_1  \|_2^2
+ t
\left( - 
\| {\varphi}_1  \|_2^2
-
\| {\varphi}_2  \|_2^2
+
2 {\varphi}_2 \cdot  {\varphi}_1
\right) 
\\[3mm]
= &
1+(t^2-t)  \| {\varphi}_2-{\varphi}_1  \|_2^2 < 1,
\\
\end{array}
$$
 for any $t\in (0,1)$. Thesis follows. 
 \qed
 \end{proof}

\subsection{Spectral  maps in annular domains}
\label{sec:affine_maps_annulus}

We denote by $\Omega = \mathcal{B}_R({0}) \setminus \mathcal{B}_r({0})$ an annular domain centered in ${0}$ with $0< r < R$ and we set $\widehat{\Omega}_{\rm pol} = (0,1)\times (-1/2,1/2)$. 
We denote by $\mathbb{P}_{J_{\rm r}} = \mathbb{P}_{J_{\rm r}}(0,1)$ the space of polynomials of degree lower or equal to $J_{\rm r}$, and by $\mathbb{F}_{J_{\rm  f}} = \mathbb{F}_{J_{\rm  f}}(-1/2,1/2)$
the Fourier space 
\begin{equation}
\label{eq:fourier_space}
\mathbb{F}_{J_{\rm  f}}  =
{\rm span}
\left\{
1, \, \cos( 2\pi x), \, \ldots, \cos( 2\pi J_{\rm  f} x),
\sin( 2\pi x), \, \ldots, \sin( 2\pi J_{\rm  f} x)
\right\}.
\end{equation}
We define 
${\Psi} : \widehat{\Omega}_{\rm pol}  \to \Omega$ and
${\Lambda} : {\Omega}  \to  \widehat{\Omega}_{\rm pol} $ such that
\begin{equation}
\label{eq:polar_map}
{\Psi}({x} = [\rho, \theta])
\,  = \,
(r + (R-r) \rho)
\left[
\begin{array}{l}
\cos (2 \pi \theta) \\
\sin (2 \pi \theta) \\
\end{array}
\right],
\quad
{\Lambda} = {\Psi}^{-1}.
\end{equation}
Then, we consider mappings of the form
\eqref{eq:spectral_map_annulus}, such that the image of $ W_M$, 
 $\mathcal{W}_M =  W_M(\mathbb{R}^M)$,  is a subset of 
\begin{subequations}
 \label{eq:ambient_space_polar}
 \begin{equation}
 \mathcal{W}_{\rm hf}^{\rm pol} = \left\{
 {\varphi} = 
 \varphi_{\rm r} \mathbf{e}_1 \, + \,
  \varphi_{\rm \theta} \mathbf{e}_2 \, : \,
   \varphi_{\rm r} \in  \mathcal{W}_{\rm hf}^{\rm pol,r} , \;\;
     \varphi_{\rm \theta} \in  \mathcal{W}_{\rm hf}^{\rm pol,\theta}
 \right\},
 \end{equation}
 with
\begin{equation}
 \mathcal{W}_{\rm hf}^{\rm pol,\theta}
\, = \, {\rm span}
\big\{
 p({x}) = p_1(x_1) p_2(x_2) \, : \, 
 p_1 \in \mathbb{P}_{J_{\rm r}}, \;
p_2 \in \mathbb{F}_{J_{\rm f}}  \big\},
\end{equation}
 and
\begin{equation}
 \mathcal{W}_{\rm hf}^{\rm pol,r}
 = 
  \big\{
 p \in   \mathcal{W}_{\rm hf}^{\rm pol,\theta} : \, 
 p(x_1, x_2 ) = 0 \; \forall \, x_1\in \{ 0, 1\}
  \big\}.
 \end{equation}
\end{subequations}

Next Proposition motivates the previous definitions.
\begin{proposition}
\label{th:annular}
Let $\Omega = \mathcal{B}_R({0}) \setminus \mathcal{B}_r({0})$. Consider the mapping  ${\Phi}$ in  \eqref{eq:spectral_map_annulus} for some $\mathbf{a} \in \mathbb{R}^M$ and 
with $\mathcal{W}_M = W_M(\mathbb{R}^M) \subset \mathcal{W}_{\rm hf}^{\rm pol}$,
(cf. \eqref{eq:ambient_space_polar}). Then, if $g_{\Phi}^{\rm pol} = {\rm det} ( \nabla {\Phi}_{\rm pol} )$ is strictly positive in $\overline{\widehat{\Omega}}_{\rm pol}$,  ${\Phi}$ is bijective in $\Omega$.
\end{proposition}

\begin{proof}
It suffices to check the hypotheses (i)-(iii) of
\cite[Proposition 2.1]{taddei2020registration}.
(i) since ${\Phi}_{\rm pol}$ is periodic in the second direction, it is easy to verify that ${\Phi}$ is smooth in $\mathcal{B}_{R+\delta}({0}) \setminus \mathcal{B}_{r-  \delta}({0})$ for $\delta<r$.
(ii) local bijectivity follows by applying the chain rule.
(iii) let ${\varphi}_{\rm pol} = W_M \mathbf{a}$; 
recalling the definition of $\mathcal{W}_{\rm hf}$ in \eqref{eq:ambient_space_polar}, we find 
${\varphi}_{\rm pol} \cdot {n} = 0$ on top and bottom edges of the unit square; the latter implies that ${\rm dist} ({\Phi}({x}), \partial \Omega  ) = 0$ for all ${x} \in \partial \Omega$.
 \qed
\end{proof}

Some comments are in order.
\begin{itemize}
\item
Proposition \ref{th:annular} does not apply to  circular domains (i.e., $r=0$) since ${\Psi}$ is singular in the origin. Nevertheless, in our numerical experience, we observe that we can construct mappings of the form \eqref{eq:spectral_map_annulus} by enforcing bijectivity at the discrete level. We refer to a future work for a thorough discussion on this issue.
\item
To deal with curved boundaries,   we consider transformations ${\Psi}, {\Lambda}$ in \eqref{eq:spectral_map_annulus} to allow parametric tangential  deformations of points on boundaries. As a result, the resulting mapping ${\Phi}$ is no longer a linear function of the mapping coefficients $\mathbf{a}$: this implies that linear reduction methods for the construction  of the reduced mapping space should be applied to ${\Phi}_{\rm pol}$.
\end{itemize}

\subsection{Spectral element maps for general domains}
\label{sec:piecewise_maps}

\begin{figure}[h!]
\centering
\includegraphics[width=0.5\textwidth]
 {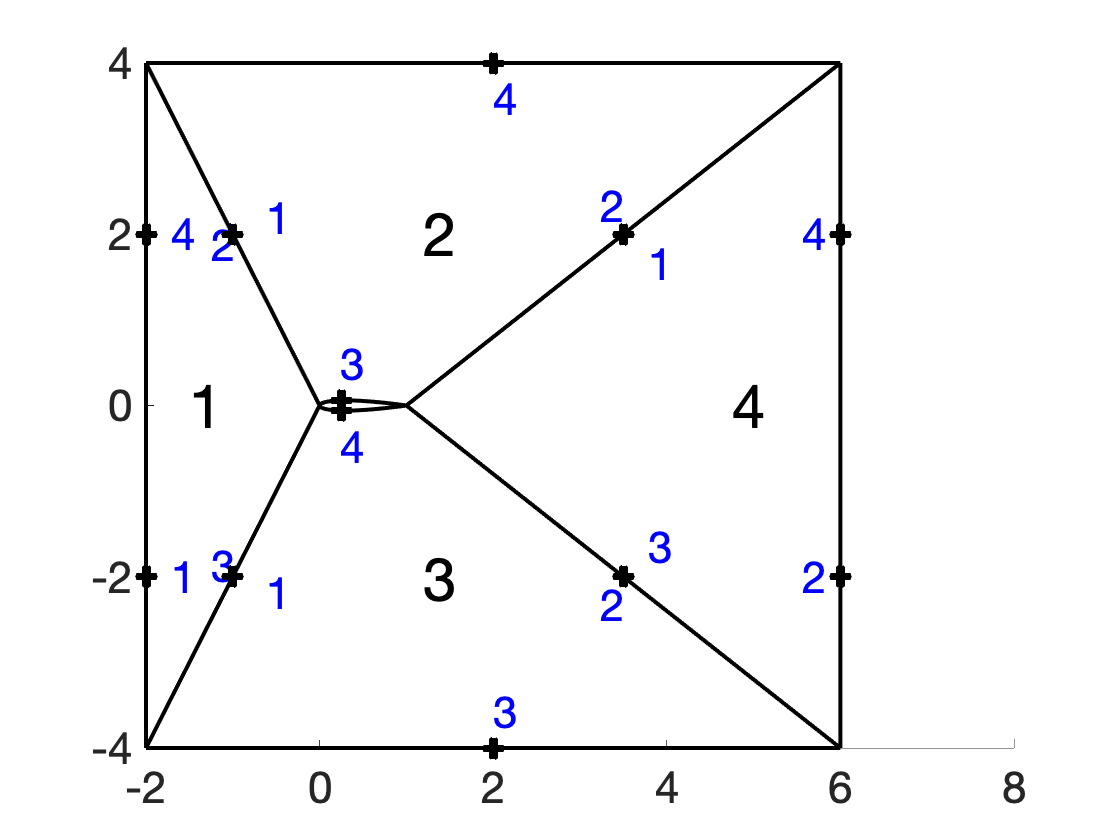}
  
 \caption{partition considered for the model problems in section \ref{sec:airfoil}  {and section \ref{sec:euler_equation}}.}
 \label{fig:vis_partitioned}
  \end{figure}

Given the domain $\Omega$, we introduce the partition $\{ \Omega_q \}_{q=1}^{N_{\rm dd}}$ and $\widehat{\Omega} = (0,1)^2$. Then, 
for $q=1,\ldots,N_{\rm dd}$, 
we define   the bijective mapping between $\widehat{\Omega}  $ and the $q$-th element of the partition
${\Psi}_q: \widehat{\Omega}  \to  \Omega_q$,
 and we denote by
${\Lambda}_q:     \Omega_q \to \widehat{\Omega}$ the inverse of ${\Psi}_q$. 
Given the tensorized polynomial space
$\mathbb{Q}_{J}$,  we introduce  the points 
$\{ \widehat{{x}}_{i,j} = (x_i^{\rm gl}, x_{j'}^{\rm gl})  \}_{i,j'=1}^J$ 
in  $\widehat{\Omega}$, where $x_1^{\rm gl},\ldots, x_{J+1}^{\rm gl}$ are the Gauss-Lobatto points  in $[0,1]$, and the corresponding Lagrangian basis
  $\{  \ell_{i,j} \}_{i,j}$ of $\mathbb{Q}_{J}$.  To fix ideas, Figure \ref{fig:vis_partitioned} shows the partition of the domain 
$\Omega$ introduced in section \ref{sec:airfoil} : the black numbers indicate the indices of the elements
$\{  \Omega_q \}_q$, while the blue numbers indicate the indices of the facets,  $\{ \partial \Omega_{q, \ell} \}_{q,\ell}$.
  
We shall consider mappings of the form \eqref{eq:piecewise_affine_mappings}:
$$
{\Phi} =
\sum_{q=1}^{N_{\rm dd}} \; \left(
{\Psi}_q \circ  {{\Phi}}_q \circ {\Lambda}_q
\right)
\mathbbm{1}_{\Omega_q},
\quad
{\rm where} \;  \;
 {{\Phi}}_q= \texttt{id} + 
 {{\varphi}}_q, \;  \;
{{\varphi}}_q \in [  \mathbb{Q}_J]^2.
$$
In order for  ${\Phi}$ to be  a mapping from $\Omega$ in itself, we should enforce the following two conditions.
\begin{enumerate}
\item
\emph{Local bijectivity:}
${\Phi}(\Omega_q) = \Omega_q$,
${\rm det} (\nabla  {\Phi}) > 0$ in $\Omega_q$ for $q=1,\ldots, N_{\rm dd}$.
\item
\emph{Continuity at interfaces:}
${\Phi} \in C(\Omega; \mathbb{R}^2)$.
\end{enumerate} 
Note that local bijectivity implies global bijectivity, but it is a much stronger condition: we here thus trade approximation power with simplicity of implementation.  We expect, however, that for a proper choice of the coarse-grained partition this family of mappings might suffice to successfully ``register" (i.e., improve the linear reducibility of) the solution manifold of interest. In the remainder of this section, we discuss how to characterize  displacements ${\varphi}_1,\ldots, {\varphi}_{N_{\rm dd}}$ that satisfy local bijectivity and global continuity.

With this in mind, we define the set  of displacements
$\overrightarrow{{\varphi}} : = [{\varphi}_1,\ldots, {{\varphi}}_{N_{\rm dd}}]$ and the space 
\begin{equation}
\label{eq:search_space_DD}
\mathcal{W}_{\rm hf}^{\rm dd} = \{
\overrightarrow{{\varphi}}   = [{\varphi}_1,\ldots, {{\varphi}}_{N_{\rm dd}}]: \;
{\varphi}_q  \in [\mathbb{Q}_{J}]^2, \;\;
q=1,\ldots,N_{\rm dd} \}.
\end{equation}
We observe that any element $\overrightarrow{{\varphi}}
= [{\varphi}_1,\ldots, {{\varphi}}_{N_{\rm dd}}]
$ in $\mathcal{W}_{\rm hf}^{\rm dd}$ is uniquely identified by the vector of coefficients $
\overrightarrow{\boldsymbol{\varphi}}
\in \mathbb{R}^{M_{\rm hf}}$, $M_{\rm hf} = 2 N_{\rm dd} (J+1)^2$, such that
\begin{equation}
\label{eq:broken_DD_algebraic}
{\varphi}_q({{X}}) 
=
\sum_{d=1}^2 \;
 \sum_{i,j=1}^{J+1}
\left( 
\overrightarrow{\boldsymbol{\varphi}}
 \right)_{\texttt{I}_{i,j,q,d}} \ell_{i,j} 
({{X}})     \; {e}_d,
\quad
{\rm for} \;  \; q=1,\ldots,N_{\rm dd},
\end{equation}
with $\texttt{I}_{i,j,q,d} = i + (j-1) (J+1) +  (J+1)^2 (q-1) +  (J+1)^2 N_{\rm dd} (d-1).$
Furthermore, 
we define the parameterizations
$\{ \widehat{{\gamma}}_{\ell} \}_{\ell=1}^4$ of the facets of $\widehat{\Omega}$, 
$ \widehat{{\gamma}}_{\ell}:[0,1] \to \partial \widehat{\Omega}$.

Recalling Proposition \ref{th:application_unit_square}, local displacement fields should satisfy
\begin{equation}
\label{eq:local_bijectivity_a}
{{\varphi}}_q \cdot \widehat{n} \big|_{\partial \widehat{\Omega}} = 0,
\quad
q=1,\ldots,N_{\rm dd},
\end{equation}
where $\widehat{n}$ denotes the outward normal to $\partial \widehat{\Omega}$. 
Furthermore, local mappings should be locally invertible: we discuss the enforcement of this condition in section \ref{sec:bijectivity}.

In order to enforce continuity at elements' interfaces, 
given 
$q \in \{1, \ldots,N_{\rm dd} \}$,
$\ell \in \{1, \ldots,4 \}$,  we denote by 
$\texttt{qext}_{\ell,q}$  the index of the neighbor element,   by 
$\texttt{ell{\_}ext}_{\ell,q}$ the index of the corresponding facet, and by 
$\texttt{orif}_{\ell,q}$ a boolean that is equal to one if the facets have the same orientation and zero otherwise. For the partition of Figure \ref{fig:vis_partitioned}, we have
$$
\begin{array}{l}
\displaystyle{
\texttt{qext} = 
\left[
\begin{array}{cccc}
-1   &   1  &   1  &   2\\
2    &  4    &  4   &  -1\\
     3  &  3   &  -1   &  3\\
-1  &   -1  &   2   &   -1\\
\end{array}
\right],
\quad
\texttt{ell{\_}ext}= 
\left[
\begin{array}{cccc}
 -1  &   2  &   3 &    2\\
 1  &  1  &   3  &  -1\\
 1   & 4  &  -1  &  2\\
-1   & -1    & 3  &   -1\\
\end{array}
\right],
} \\[6mm]
\displaystyle{
\texttt{orif}= 
\left[
\begin{array}{cccc}
\,1 \,  &  \,  1 \,&  \,  1\,&   \,  1\,\\
\, 1\, &  \,   1  \,&  \, 0\,&  \,   1\,\\
\,1 \,&    \, 1 \,  & \,  1 \,& \,   0\,\\
\,1  \,&   \,1  \, & \,  1 \,&  \,   1\,\\
\end{array}
\right].
}
\end{array}
$$
 In conclusion, if $\partial \Omega_{q, \ell} = \partial \Omega_{q', \ell'}$, $q'=\texttt{qext}_{\ell,q}$ and
 $\ell' = \texttt{ell{\_}ext}_{\ell,q}$, we obtain the conditions:
\begin{equation}
 \label{eq:continuity_interfaces}
\left\{
\begin{array}{ll}
\displaystyle{
{\varphi}_q(  \widehat{{\gamma}}_{\ell}(t)    ) =
{\varphi}_{q'}(  \widehat{{\gamma}}_{\ell'}(t)    )
}
&
\forall \; t \in [0,1], \;\;
{\rm if} \;\;\texttt{orif}_{\ell,q} = 1;
\\[3mm]
\displaystyle{
{\varphi}_q(  \widehat{{\gamma}}_{\ell}(t)    ) =
- {\varphi}_{q'}(  \widehat{{\gamma}}_{\ell'}(1-t)    )
}
&
\forall \; t \in [0,1], \;\;
{\rm if} \;\;\texttt{orif}_{\ell,q} = 0.
\\
\end{array}
\right. 
 \end{equation}
Since ${\varphi}_q, {\varphi}_{q'}$ are polynomials, it is sufficient to enforce \eqref{eq:continuity_interfaces} in the Gauss-Lobatto points.

It is interesting to characterize  the set of elements
$\mathcal{W}_{\rm hf,0}^{\rm dd}$ 
 in $\mathcal{W}_{\rm hf}^{\rm dd}$ that satisfy \eqref{eq:local_bijectivity_a} and \eqref{eq:continuity_interfaces},
\begin{equation}
\label{eq:calWhf0}
\mathcal{W}_{\rm hf,0}^{\rm dd} : =
\left\{
\overrightarrow{{\varphi}}   = [{\varphi}_1,\ldots, {{\varphi}}_{N_{\rm dd}}]: \;
\overrightarrow{{\varphi}}   \;\;  {\rm satisfies} \;\; 
\eqref{eq:local_bijectivity_a}  - \eqref{eq:continuity_interfaces}
\right\}.
\end{equation} 
   Conditions
\eqref{eq:local_bijectivity_a} correspond to $4 (J+1) N_{\rm dd}$ linear conditions for $\overrightarrow{\boldsymbol{\varphi}}$, while \eqref{eq:continuity_interfaces} corresponds to $N_{\rm int} (J-1)$ additional conditions, where $N_{\rm int}$ denotes the number of interior facets. In conclusion, we find that 
$\mathcal{W}_{\rm hf,0}$   is a linear space of dimension
$(2 (J+1)^2 - 4 (J+1) ) N_{\rm dd} -  (J-1) N_{\rm int}$. We might represent elements of $\mathcal{W}_{\rm hf,0}^{\rm dd}$ as the sum of 
$N_{\rm dd}$ local displacements that are clamped at the boundary of $\widehat{\Omega}$, and of
$N_{\rm f}$ tangential displacements, where $N_{\rm f}$ denotes the number of facets of the partition.

\subsection{Practical enforcement of the bijectivity condition}
\label{sec:bijectivity}

Given the spectral or spectral element mapping ${\Phi} = {\Phi}(\cdot; \mathbf{a})$, in view of registration, we shall devise conditions for the mapping coefficients $\mathbf{a}$ to ensure that 
(i) ${\rm det} \nabla {\Phi}(\cdot ; \mathbf{a}) > 0$ in $\overline{\Omega}$,  and
(ii) given the mesh $\mathcal{T}_{\rm hf}$, the mapped mesh ${\Phi}(\mathcal{T}_{\rm hf}, \mathbf{a})$ is   well-defined  (cf. Definition \ref{def:discrete_bijectivity}).  Note that (i) guarantees that ${\Phi}$ is bijective in $\Omega$, while (ii) guarantees that  
${\Phi}(\mathcal{T}_{\rm hf}, \mathbf{a})$ can be used for FE calculations.
We restrict our attention to the general case of two-dimensional domains; the other two cases can be treated similarly.

In order to enforce (i), 
since $\{  {\Psi}_q \}_{q=1}^{N_{\rm dd}}$ are bijective, it suffices to enforce that the internal maps ${{\Phi}}_1,\ldots, {{\Phi}}_{N_{\rm dd}}$ are invertible. Following \cite{taddei2020registration,taddei2021space}, we introduce 
\begin{subequations}
\label{eq:bijectivity_calC}
\begin{equation}
\mathfrak{C}(  \mathbf{a} ) \; : = \;
\sum_{q=1}^{N_{\rm dd}} \; 
\int_{\widehat{\Omega}} \;
{\rm exp} \left(
\frac{\epsilon  -   \widehat{g}_{\Phi_q}({x}) }{  C_{\rm exp}  }
\right)
\,+ \,
{\rm exp} \left(
\frac{ \widehat{g}_{\Phi_q}({x}) - 1/\epsilon   }{  C_{\rm exp}  }
\right)
\; d{x}
-\delta N_{\rm dd};
\end{equation}
where $\widehat{g}_{\Phi_q}  : = {\rm det} \nabla {\Phi}_q$ for $q=1,\ldots,N_{\rm dd}$. We choose $C_{\rm exp}, \epsilon,\delta$ as follows:
\begin{equation}
\epsilon=0.1,  \quad
C_{\rm exp} = 0.025 \epsilon, \quad
\delta = 1.
\end{equation}
\end{subequations}
As explained in \cite[section 2.2]{taddei2020registration}, if $\mathfrak{C}(  \mathbf{a} )$ is negative, the local maps $\{  {\Phi}_q(\cdot ; \mathbf{a}) \}_q$ are bijective, provided that $\{ \|  \nabla \widehat{g}_{\Phi_q} \|_{L^{\infty}(\widehat{\Omega})  } \}_q$ are all moderate.

To ensure discrete bijectivity, we introduce the mesh distortion indicator
\begin{subequations}
\begin{equation}
\label{eq:local_mesh_distortion}
\mathfrak{f}_{\rm msh,k}
: = 
\frac{1}{2}
\frac{\|  \nabla {\Psi}_{k,\Phi}^{\rm hf,1}   \|_{\rm F}^2 }{
|  {\rm det} (  \nabla {\Psi}_{k,\Phi}^{\rm hf,1}   )  |},
\;k=1,\ldots,N_{\rm e},
\end{equation}
where 
$\| \cdot \|_{\rm F}$ is the Frobenius norm and 
${\Psi}_{k,\Phi}^{\rm hf,1} $ is the elemental mapping \eqref{eq:psi_mapping_mapped} associated with  a \texttt{p}=1 discretization. We observe that the indicator  \eqref{eq:local_mesh_distortion} is widely used for high-order mesh generation, and has also been considered in  \cite{zahr2020implicit} to prevent mesh degradation, in the DG framework. Then, we introduce the integral function:
\begin{equation}
\label{eq:registration_statement_e}
\mathfrak{R}_{\rm msh} (\mathbf{a} ) 
=
\sum_{k=1}^{N_{\rm e}} \; |\texttt{D}_k |
{\rm exp} \left(
\mathfrak{f}_{\rm msh,k}\left( 
{\Phi}(\cdot; \mathbf{a})
\right)
\,-\,
\mathfrak{f}_{\rm msh,max}
\right).
\end{equation}
Here, $\mathfrak{f}_{\rm msh,max}>0$ is a given threshold;
in all our numerical experiments, we set $\mathfrak{f}_{\rm msh,max}=10$.
\end{subequations}

\subsection{Implementation considerations and extension to parameterized geometries}
 \label{sec:implementation_param}

As explained in section \ref{sec:approx_form_gen}, it is of paramount importance to rapidly deform the mesh points $\{ {x}_j^{\rm hf}  \}_j$ for a new value of the parameters $\mu \in \mathcal{P}$. Towards this end, following \cite{taddei2020discretize}, given the mesh $\mathcal{T}_{\rm hf}$ and the partition $\{ \Omega_q  \}_{q=1}^{N_{\rm dd}}$, we compute 
$\texttt{I}_{\Phi} \in \{ 1,\ldots,N_{\rm dd} \}^{N_{\rm hf}}$ 
and $\{  {x}_j^{\rm hf,ref} \}_{j=1}^{N_{\rm hf}} \subset [0,1]^2$ such that
$( \texttt{I}_{\Phi}  )_j$ denotes the label of the region to which the $j$-th node of the mesh ${x}_j^{\rm hf}$ belongs, and 
${x}_j^{\rm hf,ref} :=
{\Psi}_{( \texttt{I}_{\Phi}  )_j    }^{-1} \left( {x}_j^{\rm hf}  \right)$.
Then, given $\mathbf{a} \in \mathbb{R}^M$, we define the deformed nodes using the identity:
\begin{equation} 
\label{eq:onlinePhi}
{\Phi}_{\mu}( {x}_j^{\rm hf}     )
=
{\Psi}_{ ( \texttt{I}_{\Phi}  )_j    } 
\left( 
{x}_j^{\rm hf,ref} 
+ 
W_M(  {x}_j^{\rm hf,ref}  ) \mathbf{a}
 \right),
\quad
j=1,\ldots,N_{\rm hf}.
\end{equation}
The latter can be rapidly evaluated, provided that ${\Psi}_1,\ldots,{\Psi}_{N_{\rm dd}}$ can be computed in $\mathcal{O}(1)$ flops.

We also observe that \eqref{eq:piecewise_affine_mappings}  can be trivially modified to deal with parameterized geometries. We define the parametric partition $\{ \Omega_{q,\mu}  \}_{q=1}^{N_{\rm dd}}$ of $\Omega_{\mu}$ and the maps ${\Psi}_q: \widehat{\Omega} \times \mathcal{P} \to \Omega_{q,\mu}$ and their inverses ${\Lambda}_q:  \Omega_{q,\mu} \times \mathcal{P} \to \widehat{\Omega}$, for 
$q=1,\ldots,N_{\rm dd}$. Then, given $\mu, \bar{\mu} \in \mathcal{P}$, we define the bijection 
${\Phi}$ from $\Omega_{\bar{\mu}}$ to $\Omega_{\mu}$ such that
\begin{equation}
\label{eq:partitioned_approach_paramgeo}
{\Phi} =
\sum_{q=1}^{N_{\rm dd}} \; \left(
{\Psi}_{q,\mu} \circ  {{\Phi}}_q \circ {\Lambda}_{q,\bar{\mu}}
\right)
\mathbbm{1}_{\Omega_{q,\bar{\mu}}},
\end{equation}
where  $ {{\Phi}}_q= \texttt{id} + W_M^q \mathbf{a}: \widehat{\Omega}\to \widehat{\Omega}$. Note that the inverses $\{  {\Lambda}_{q,\bar{\mu}}  \}_q$ are computed for a select value of the parameter; we can thus modify
\eqref{eq:onlinePhi}  to rapidly deform the mesh $\mathcal{T}_{\rm hf}$ for any  given set of mapping coefficients $\mathbf{a}$ and any  $\mu \in \mathcal{P}$.

{The use of Gordon-Hall maps relies on the assumption that explicit parameterizations of the boundary are  available.  In the model reduction literature, several authors have proposed different strategies to deal with more complex domains, particularly for vascular applications
(e.g., \cite{lassila2014model,manzoni2012model}). Our approach does not require the use of Gordon-Hall maps; we might indeed build $\{ \Psi_q  \}_{q=1}^{N_{\rm dd}}$ using other geometry registration techniques that 
are better suited for the particular geometry of interest.
Note, however, that continuity at interfaces (cf.\eqref{eq:continuity_interfaces}) requires the compatibility of maps of neighboring elements at the shared interface:  enforcement of this condition is trivial for Gordon-Hall maps, while it might be more involved for other geometry  registration techniques. 
}

\section{Registration}
\label{sec:registration}
We  adapt the registration strategy proposed in \cite{taddei2021space} to the more general framework considered in this paper. We here discuss the approach  for general two-dimensional domains (cf. section \ref{sec:piecewise_maps}): the other two cases can be handled similarly.  We consider the case of parameterized geometries: given the family of parameterized domains  $\{ \Omega_{\mu} : \mu \in \mathcal{P}  \}$ and the parametric partition
$\{ \Omega_{q, \mu} \}_{q=1}^{N_{\rm dd}}$,  $\bar{\mu} \in \mathcal{P}$, 
we set  $\Omega := \Omega_{\bar{\mu}}$ and $\{ \Omega_{q}:= \Omega_{q, \bar{\mu}} \}_{q=1}^{N_{\rm dd}}$; furthermore, we denote by $\mathcal{T}_{\rm hf}$ the FE mesh in $\Omega$.

\subsection{Parametric registration}
\label{sec:parametric_reg}
Exploiting the discussion in section \ref{sec:piecewise_maps}, we can rewrite \eqref{eq:piecewise_affine_mappings} as
\begin{equation}
\label{eq:calN_expansion}
{\mathcal{N}}(\cdot;  \, \overrightarrow{{\varphi}}) \, = \,
\sum_{q=1}^{N_{\rm dd}} \; \left(
{\Psi}_{q,\mu} \circ  {{\Phi}}_q \circ {\Lambda}_{q}
\right) \mathbbm{1}_{\Omega_q},
\quad
{{\Phi}}_q = \texttt{id} + {\varphi}_q,
\;\;
\overrightarrow{{\varphi}} \in \mathcal{W}_{\rm hf,0}^{\rm dd},
\end{equation}
where $\{ {\Lambda}_q := {\Lambda}_{q,\bar{\mu}} \}$.
 In view of dimensionality reduction, we equip $\mathcal{W}_{\rm hf,0}^{\rm dd}$ with the norm
 \begin{equation}
 \label{eq:norm_space}
\vertiii{  \overrightarrow{{\varphi}}   }^2  : =((\overrightarrow{{\varphi}} , \overrightarrow{{\varphi}}  ))
=
\sum_{q=1}^{N_{\rm dd}}  | \Omega_q  |
 \|  {\varphi}_q  \|_{H^2( \widehat{\Omega}  )}^2.
 \end{equation}
Furthermore, given the space $\mathcal{W}_M \subset \mathcal{W}_{\rm hf,0}^{\rm dd}$, $M = {\rm dim} ( \mathcal{W}_M  ) \leq  {\rm dim} ( \mathcal{W}_{\rm hf,0}^{\rm dd} ) = M_{\rm hf} $, we introduce the isometry $W_M: \mathbb{R}^M \to \mathcal{W}_M$, $\vertiii{W_M \mathbf{a} } = \|  \mathbf{a} \|_2$ for all $\mathbf{a} \in \mathbb{R}^M$.

We define the registration sensor ${s}: \mu \in \mathcal{P} \to [ L^2(\widehat{\Omega})  ]^{N_{\rm dd}}$ and the manifold $\mathcal{M}_{\rm s} = \{ {s}_{\mu}: \mu \in \mathcal{P}   \}$. The field ${s}_{\mu}$ might be an explicit function of the solution field $u_{\mu}$ or of certain parametric coefficients associated with the PDE.
As  discussed  in \cite[Remark 3.1]{taddei2021space}, the sensor should capture relevant features associated with the solution field.

The first step of the registration procedure is to devise an automatic procedure  to learn an optimal displacement $\overrightarrow{{\varphi}}^{\star}$ based on a target ${s}^{\star} \in  [ L^2(\widehat{\Omega})  ]^{N_{\rm dd}}$, an $N$-dimensional approximation space 
$\mathcal{S}_N \subset [ L^2(\widehat{\Omega})  ]^{N_{\rm dd}}$ (referred to as \emph{template space}), an $M$-dimensional approximation space  $\mathcal{W}_M$ for the mapping, and a mesh $\mathcal{T}_{\rm hf}$ of the domain $\Omega$:
\begin{equation}
\label{eq:registration_function}
\left[ \overrightarrow{{\varphi}}^{\star}, \mathfrak{f}_{N,M}^{\star}  \right]
\, = \,
\texttt{registration} \left(
{s}^{\star}, \mathcal{S}_N, \mathcal{W}_M, \mathcal{T}_{\rm hf}
\right).
\end{equation}
The displacement field $\overrightarrow{{\varphi}}^{\star}$ should guarantee discrete bijectivity with respect to $\mathcal{T}_{\rm hf}$ and also that  ${\mathcal{N}}(\cdot;  \, \overrightarrow{{\varphi}})$ is bijective in $\Omega$. 
Here, $\mathfrak{f}_{N,M}^{\star} >0$ measures performance of registration.
We discuss in detail the registration algorithm \eqref{eq:registration_function}  in section \ref{sec:optimization}.

Given \eqref{eq:registration_function} and the set of snapshots $\mathcal{M}_{\rm s}$, $\{ {s}^k ={s}_{\mu^k}   \}_{k=1}^{n_{\rm train}}$, we resort to the greedy procedure in \cite[Algorithm 1]{taddei2021space} to simultaneously build the spaces $\mathcal{S}_N$, $\mathcal{W}_M$ and the mapping coefficients $\{  \mathbf{a}^k \}_k$ associated with  
the orthonormal  basis   $\{ \overrightarrow{{\varphi}}_m  \}_{m=1}^M$ of $ \mathcal{W}_M$ --- 
$\overrightarrow{{\varphi}}^{k, \star} = 
\sum_m (  \mathbf{a}^k  )_m\overrightarrow{{\varphi}}_m
=  W_M \mathbf{a}^k$. For completeness, we summarize the parametric registration algorithm in Algorithm \ref{alg:registration}. The function
$$
[  W_M , 
\; 
\{  \mathbf{a}^k  \}_k ]  =
\texttt{POD} \left( \{ \overrightarrow{{\varphi}}^{k, \star}  \}_{k=1}^{n_{\rm train}}, 
tol_{\rm pod}  , (( \cdot, \cdot ))   \right)
$$
corresponds to the application of POD based on the norm \eqref{eq:norm_space}; given the eigenvalues $\{ \lambda_m \}_m$ of the Gramian matrix and the tolerance $tol_{\rm pod} > 0$, the size $M$ of the space is built based on the criterion
\begin{equation}
\label{eq:POD_cardinality_selection}
M := \min \left\{
M': \, \sum_{m=1}^{M'} \lambda_m \geq  \left(1 - tol_{\rm pod} \right) 
\sum_{i=1}^{n_{\rm train}} \lambda_i
\right\}.
\end{equation} 
The coefficients 
 $\{ \mathbf{a}^k \}_k$ are given by
$ ( \mathbf{a}^k   )_m = ((   W_M \mathbf{e}_m ,  \overrightarrow{{\varphi}}^{k, \star} ))$ for $m=1,\ldots,M$ and $k=1,\ldots,n_{\rm train}$, where $\mathbf{e}_1,\ldots,\mathbf{e}_M$ are the vectors of the canonical basis in $\mathbb{R}^M$.

Given the dataset $\{ (\mu^k, \mathbf{a}^k)  \}_{k=1}^{n_{\rm train}}$, 
following \cite{taddei2020registration,taddei2021space},  
we resort to a multi-target regression algorithm to learn a regressor $\mu \mapsto \widehat{\mathbf{a}}_{\mu}$ and ultimately the  parametric mapping
\begin{equation}
\label{eq:parametric_mapping_Phi}
{\Phi} : \Omega \times \mathcal{P}  \to \mathbb{R}^2,
\quad
{\Phi}_{\mu} : =
{\mathcal{N}}(\cdot ; W_M \widehat{\mathbf{a}}_{\mu} ).
\end{equation}
We here resort to radial basis function (RBF, \cite{wendland2004scattered}) approximation: other regression algorithms could also be considered. We observe that purely data-driven regression techniques  do not enforce bijectivity for out-of-sample parameters: in practice, we should thus consider sufficiently large training sets in Algorithm \ref{alg:registration}.

\begin{algorithm}[H]                      
\caption{Registration algorithm}     
\label{alg:registration}     

 \small
\begin{flushleft}
\emph{Inputs:}  $\{ (\mu^k,  {s}_{\mu^k}) \}_{k=1}^{n_{\rm train}} \subset  \mathcal{P} \times \mathcal{M}_{\rm s}$ snapshot set, 
$\mathcal{S}_{N_0} = {\rm span} \{ {\psi}_n \}_{n=1}^{N_0}$ initial template space;
\smallskip


\emph{Outputs:} 
${\mathcal{S}}_N = {\rm span} \{ {\psi}_n  \}_{n=1}^N$ template space, 
$\mathcal{W}_M = {\rm span} \{ \overrightarrow{{\varphi}}_m  \}_{m=1}^M$ displacement space,
$\{  \mathbf{a}^k  \}_k$ mapping coefficients.
\end{flushleft}                      

 \normalsize 

\begin{algorithmic}[1]
\State
Set
$\mathcal{S}_{N=N_0} = \mathcal{S}_{N_0}$,
$\mathcal{W}_M =\mathcal{W}_{\rm hf,0}^{\rm dd}$.
\vspace{3pt}

\For {$N=N_0, \ldots, N_{\rm max}-1$ }

\State
$
\left[ \overrightarrow{{\varphi}}^{\star,k}, \mathfrak{f}_{N,M}^{\star,k}  \right]
\, = \,
\texttt{registration} \left(
{s}^{k}, \mathcal{S}_N, \mathcal{W}_M, \mathcal{T}_{\rm hf}
\right)$ for $k=1,\ldots,n_{\rm train}$.
\vspace{3pt}

\State
$[  W_M , \; \{  \mathbf{a}^k  \}_k ]  =
\texttt{POD} \left( \{ \overrightarrow{\varphi}^{\star,k} \}_{k=1}^{n_{\rm train}}, 
tol_{\rm pod}  , \vertiii{\cdot}   \right)$, $\mathcal{W}_M = W_M(\mathbb{R}^M)$.
\vspace{3pt}

\If{   $\max_k   \mathfrak{f}_{N,M}^{\star,k}   < \texttt{tol}$}, \texttt{break}

\Else

\State
$\mathcal{S}_{N+1} = \mathcal{S}_{N} \cup {\rm span} 
\{ {s}_{\mu^{k^{\star}}}    \circ  {\Phi}^{\star,k^{\star}}  \}$ with 
$k^{\star} = {\rm arg} \max_k \mathfrak{f}_{N,M}^{\star,k}$.
\EndIf
\EndFor
\end{algorithmic}
\end{algorithm}

After having built the parametric mapping ${\Phi}$, we might resort to standard (linear) reduction techniques to compute the linear approximation $\mu \mapsto \widehat{\mathbf{u}}_{\mu} = \mathbf{Z}_N \widehat{\boldsymbol{\alpha}}_{\mu}$ in 
\eqref{eq:registration_fundamental_approx}.
As stated in the introduction, we here resort to a fully non-intrusive strategy based on POD and RBF approximation:
(i) we generate snapshots $\{  \mathbf{u}_{\mu^k}^{\rm hf}   \}_{k}$ by solving the parametric differential problem in the mapped meshes
$\mathcal{T}_{\rm hf,\mu^k} = {\Phi}_{\mu^k}(\mathcal{T}_{\rm hf}  )$, $k=1,\ldots,n_{\rm train}$;
(ii) we apply POD to build the linear operator $Z_N$ and the solution coefficients 
$\{ \boldsymbol{\alpha}^k = \mathbf{Z}_N^T \mathbf{X} \mathbf{u}_{\mu^k}^{\rm hf}  \}_{k=1}^{n_{\rm train}}$,
where the matrix $\mathbf{X}$ is defined in \eqref{eq:Xnorm};
(iii) we use RBF regression to learn the regressor $\mu \mapsto \widehat{\boldsymbol{\alpha}}_{\mu}$ based on the dataset
$\{  ( \mu^k,  \boldsymbol{\alpha}^k )  \}_k$. 
If  snapshots $\{  u_{\mu^k}^{\rm hf}   \}_{k}$ have already been computed to generate 
$\{  {s}_{\mu^k}^{\rm hf}   \}_{k}$, we might resort to mesh interpolation to generate the snapshot set $\{  \mathbf{u}_{\mu^k}^{\rm hf}   \}_{k}$.
We refer to the above-mentioned papers, for an application of projection-based pMOR techniques to approximate the mapped field.

\begin{remark}
\label{remark:norm_map}
\emph{Choice of the mapping norm.}
In our experience, for highly-anisotropic partitions,   it might be important to  consider    mapping norms that take  into account the shapes of the partition elements $\{ \Omega_q \}_q$. We propose to introduce the linear approximations
$\{  \widetilde{{\Psi}}_q \}_q$ of the mappings $\{  {\Psi}_{q,\bar{\mu}} \}_q$ and then consider the norm
\begin{equation}
\label{eq:norm_map_modified}
\vertiii{  \overrightarrow{{\varphi}}   }^2  : =((\overrightarrow{{\varphi}} , \overrightarrow{{\varphi}}  ))
=
\sum_{q=1}^{N_{\rm dd}}  | \Omega_q  |
 \| \nabla  \widetilde{{\Psi}}_q   \left(   {\varphi}_q \circ \widetilde{{\Psi}}_q^{-1}   \right)  \|_{H^2( \widetilde{\Omega}_q  )}^2,
 \quad
 \widetilde{\Omega}_q = 
 \widetilde{{\Psi}}_q(  \widehat{\Omega}  ).
\end{equation}
This choice is simple to implement and performs well in practice; we refer to a future work for a thorough assessment of the choice of the mapping norm.
\end{remark}

\subsection{Optimization statement}
\label{sec:optimization}

Given the target ${s} = [s_1,\ldots,s_{N_{\rm dd}}] \in [ L^2(\widehat{\Omega}) ]^{N_{\rm dd}}$, the parameter $\mu \in \mathcal{P}$, the template space $\mathcal{S}_N$, the mapping space $\mathcal{W}_M = W_M(\mathbb{R}^M)$ and the mesh $\mathcal{T}_{\rm hf}$, we define $\overrightarrow{\varphi}^{\star} = W_M \mathbf{a}^{\star}$, as the solution to the optimization problem
\begin{subequations}
\label{eq:optimization_statement}
\begin{equation}
\begin{array}{l}
\displaystyle{
\min_{\mathbf{a} \in \mathbb{R}^M}
\mathfrak{f}( \mathbf{a};   {s}, \mathcal{S}_N, \mathcal{W}_M  ) 
\; + \;
\xi \| \mathbf{A}_{\rm stab}^{1/2} \mathbf{a}  \|_2^2 
\; + \;
\xi_{\rm msh} \mathfrak{R}_{\rm msh} (\mathbf{a}; \mu );
} \\[3mm]
\displaystyle{
{\rm subject \; to} \;\;
\mathfrak{C}(\mathbf{a}) \leq 0.
}
\\
\end{array}
\end{equation}
Here, $\mathfrak{f}$ measures the projection error associated with the mapped target ${s}^{\star}$ with respect to the template space $\mathcal{S}_N$,
\begin{equation}
\label{eq:proximity_measure}
\mathfrak{f}( \mathbf{a};   {s}^{\star}, \mathcal{S}_N, \mathcal{W}_M  ) \; := \;
\min_{{\psi} \in \mathcal{S}_N }
\sum_{q=1}^{N_{\rm dd}} \int_{\widehat{\Omega}} \;
\left(
s_q \circ {\Phi}_q(\cdot; \mathbf{a})
- \psi_q
\right)^2 g_{\Psi_{q,\mu}} \; d{x},
\end{equation}
where  ${\Phi}_q(\cdot; \mathbf{a}) = \left(
\texttt{id} + W_M \mathbf{a} \right)_q$, 
$q=1,\ldots,N_{\rm dd}$:
we thus set  $\mathfrak{f}_{N,M}^{\star} =
\mathfrak{f}^{\star}( \mathbf{a}^{\star};   {s}^{\star}$, $\mathcal{S}_N, \mathcal{W}_M)$ in \eqref{eq:registration_function}.
The matrix  $\mathbf{A}_{\rm stab} \in \mathbb{R}^{M,M}$ is the symmetric  positive semi-definite matrix  associated with the $H^2(\widehat{\Omega})$-seminorm:
\begin{equation}
\| \mathbf{A}_{\rm stab}^{1/2} \mathbf{a}  \|_2^2 
=
\sum_{q=1}^{N_{\rm dd}} 
\big|  (W_M \mathbf{a})_q    \big|_{H^2(\widehat{\Omega})}^2;
\end{equation} 
and  $\mathfrak{R}_{\rm msh}$, $\mathfrak{C}$ are the functionals introduced in section \ref{sec:bijectivity} to enforce discrete and continuous bijectivity.
\end{subequations}
 
 The functional  $\mathfrak{f}$ is designed to measure the approximability of the solution $u$ through a low-dimensional --- yet to be determined --- linear space in the reference configuration. To motivate this claim, consider 
 $$
 {s}_{\mu}: =
 \left[
 u_{\mu}\circ {\Psi}_{1,\mu}, \ldots,
  u_{\mu}\circ {\Psi}_{N_{\rm dd},\mu}
 \right]^T
 =
 \mathfrak{s}_{\mu}(u_{\mu}),
  $$ 
 and set $\mathcal{S}_N = {\rm span} \{ {s}_{\bar{\mu}} ( \zeta_n  ) \}_{n=1}^N  \subset \mathcal{X}$, where
$\mathcal{Z}_N = {\rm span} \{ \zeta_n \}_{n=1}^N$  
  is an approximation space for the (yet to be built) mapped solution manifold $\widetilde{\mathcal{M}}$. Then, for any $\mathbf{a} \in \mathbb{R}^M$, we find
$$
\begin{array}{l}
\displaystyle{
\min_{\zeta \in \mathcal{Z}_N} 
\int_{\Omega_{\bar{\mu}}} \;
\left(
u_{\mu} \circ {\Phi}_{\mu} - \zeta \right)^2 
d{x}
=
\sum_{q=1}^{N_{\rm dd}} 
\int_{\Omega_{q,\bar{\mu}}} \;
\left(
u_{\mu} \circ {\Psi}_{q,\mu} \circ {\Phi}_{q,\mu} \circ {\Lambda}_q - \zeta \right)^2 
d{x}
}
\\[3mm]
\displaystyle{
=
\min_{ \zeta \in \mathcal{Z}_N}
\; 
\sum_{q=1}^{N_{\rm dd}} 
 \int_{ \widehat{\Omega} }  \;
\;
\left(
s_{q,\mu}  \circ {\Phi}_{q,\mu}    - \zeta \circ \Psi_{q,\bar{\mu}}  \right)^2  g_{\Psi_{q,\bar{\mu}}} d{x}
}
\\[3mm]
\displaystyle{
=
\min_{ {\psi}    \in \mathcal{S}_N}
\; 
\sum_{q=1}^{N_{\rm dd}} 
 \int_{ \widehat{\Omega} }  \;
\;
\left(
s_{q,\mu}  \circ {\Phi}_{q,\mu}    - \psi_q \right)^2  g_{\Psi_{q,\bar{\mu}}} d{x},
}
\\
\end{array}
$$
which is $\mathfrak{f}( \mathbf{a};   {s}, \mathcal{S}_N, \mathcal{W}_M  ) $.
We discuss in detail how to choose the registration sensor in
{the next section}: here, we note that it is important to define ${s}$ over a structured grid, to reduce the cost of computing $s_1, \ldots,s_{N_{\rm dd}}$ at the mapped quadrature points.

We observe that \eqref{eq:optimization_statement} depends on several hyper-parameters. We   discussed in section \ref{sec:spectral_maps} the choices of $\epsilon, \delta, C_{\rm exp}$ in $\mathfrak{C}$,  and of $\mathfrak{f}_{\rm msh,max}$ in $\mathfrak{R}_{\rm msh}$; the choices of $\xi$ and $\xi_{\rm msh}$ balance accuracy --- measured by $\mathfrak{f}$ --- and smoothness of the mapping --- measured by $\| \mathbf{A}_{\rm stab}^{1/2}\cdot \|_2^2$ --- and of the mesh --- measured by $\mathfrak{R}_{\rm msh}$. In our experience, the  choice of $\xi$ is critical for performance and should be carefully tuned: we discuss the choice of the hyper-parameters for the considered model problems in section \ref{sec:numerics}. Since the registration problem is non-convex, careful initialization of the iterative optimization algorithm is important: we refer to  \cite[section 3.1.2]{taddei2020registration} for further details.

We also remark that \eqref{eq:optimization_statement} differs from the registration statement in \cite{taddei2021space} due to the addition of the penalty term $\mathfrak{R}_{\rm msh}$ in the objective function. As observed in section \ref{sec:spectral_maps}, this term enforces discrete bijectivity  and is particularly relevant to cope with anisotropic meshes.

\subsection{Choice of the registration sensor}
\label{sec:registration_sensor}

{We denote by $\widehat{\mathcal{T}}_{\rm hf} = \left(
\{ \widehat{x}_{j}^{\rm hf}  \}_{j=1}^{\widehat{N}_{\rm hf}}, \widehat{\texttt{T}} \right)$ a structured mesh of $\widehat{\Omega}$ and we denote by $\widehat{\mathfrak{X}}_{\rm hf}$ the continuous FE space of order \texttt{p} associated with 
$\widehat{\mathcal{T}}_{\rm hf}$. 
In the case of geometrical parameterizations, we define   the  purely-geometrical mapping
$\Phi_{\mu}^{\rm geo}  = \sum_{q=1}^{N_{\rm dd}}
\Psi_{q,\mu}\circ \Lambda_q \mathbbm{1}_{\Omega_q} : \Omega \to \Omega_{\mu}$, which corresponds to the choice $\Phi_1 = \ldots = \Phi_{N_{\rm dd}} = \texttt{id}$ in \eqref{eq:partitioned_approach_paramgeo}.
Given the solution field 
$( {\mathcal{T}}_{\rm hf,\mu}^{\rm geo} = \Phi_{\mu}^{\rm geo} (  \mathcal{T}_{\rm hf} ) , \mathbf{u}_{\mu}^{\rm hf})$ for some $\mu\in \mathcal{P}$, we  discuss in this section how to compute the registration sensor $s_{\mu}$.  To shorten notation, we assume that 
$s_{\mu}$ is computed based on $\mathbf{u}_{\mu}^{\rm hf}$: in the example of section \ref{sec:euler_equation}, we compute 
$s_{\mu}$ based on a scalar function of  $\mathbf{u}_{\mu}^{\rm hf}$, the Mach number.

We investigate two separate strategies. In the first approach, we solve a smoothing problem in each element of the partition to obtain a smooth projection of $u_{\mu}\circ \Psi_{\mu}$ on 
$\widehat{\mathfrak{X}}_{\rm hf}$; given $q=1,\ldots,N_{\rm dd}$, we define $( s_{\mu} )_q$ such that
\begin{equation}
\label{eq:sensor_approach1}
\left( s_{\mu}  \right)_q := {\rm arg} \min_{\varphi  \in \widehat{\mathfrak{X}}_{\rm hf} } \;
\xi_{\rm s} \| \nabla \varphi  \|_{L^2(\widehat{\Omega})}^2
\, + \,
\sum_{j: ( \texttt{I}_{\Phi} )_j = q    } \, 
\left(
\varphi ( {x}_j^{\rm hf, ref}  )
-
\left( \mathbf{u}_{\mu}^{\rm hf} \right)_j
\right)^2,
\end{equation}
where $\xi_{\rm s}>0$ is a smoothing parameter. In the second approach, we first preprocess  the field $\widetilde{u}_{\mu}^{\rm hf} \in \mathfrak{X}_{\mathcal{T}_{\rm hf}}$ associated with 
the FE vector $\mathbf{u}_{\mu}^{\rm hf}$
(cf. \eqref{eq:vector2field})
 by solving the smoothing problem
\begin{subequations}
\label{eq:sensor_approach2}
\begin{equation}
\int_{\Omega}  \xi_{\rm s} \nabla \widetilde{u}_{\mu}^{\rm sm} \cdot \nabla v \; + \;
(  \widetilde{u}_{\mu}^{\rm sm} - \widetilde{u}_{\mu}^{\rm hf}  ) \, v
\; d\,x = 0,
\quad \forall \, v\in  
\mathfrak{X}_{\mathcal{T}_{\rm hf}} \cap H_0^1(\Omega),
\;\;
\widetilde{u}_{\mu}^{\rm sm}|_{\partial \Omega} = \widetilde{u}_{\mu}^{\rm hf},
\end{equation}
with   $\xi_{\rm s}>0$; then, we define 
\begin{equation}
\label{eq:sensor_approach2b}
\left( s_{\mu}  \right)_q (  \widehat{x}_j^{\rm hf}   )
=
u_{\mu}^{\rm sm} \left(   \Psi_{q,\mu} (  \widehat{x}_j^{\rm hf}     )   \right)
=
\widetilde{u}_{\mu}^{\rm sm} \left(   \Psi_{q,\bar{\mu}} (  \widehat{x}_j^{\rm hf}     )   \right).
\end{equation}
Note that in the last equality we used the fact that
$$
\widetilde{u}_{\mu}^{\rm sm} \circ \Psi_{q,\bar{\mu}} 
=
{u}_{\mu}^{\rm sm} \circ \Psi_{q,\mu} \circ \Lambda_{q} 
\circ 
\Psi_{q,\bar{\mu}} 
=
{u}_{\mu}^{\rm sm} \circ \Psi_{q,\mu}.
$$
\end{subequations}
Provided that the snapshots are all computed using the meshes 
$\{ {\mathcal{T}}_{\rm hf,\mu}^{\rm geo} : \mu \in \mathcal{P} \}$, evaluation of \eqref{eq:sensor_approach2b} for all snapshots in the training set requires mesh interpolation over $N_{\rm dd} \cdot \widehat{N}_{\rm hf}$ nodes.

In our experience, the two approaches  lead to very similar results for problems with quasi-uniform meshes and smooth fields (see the examples in sections \ref{sec:heat_transfer} and \ref{sec:airfoil}); on the other hand, the first  approach might lead to excessive oscillations for meshes that are refined in specific regions of the domain (see the example in section \ref{sec:euler_equation}). For this reason, we resort to the former for the first two examples in section \ref{sec:numerics} and to the latter for the third example.
We further remark that, in the case of annular domains,
the mapping ${\Phi}_{\rm pol}$ is not guaranteed to map $\widehat{\Omega}_{\rm pol}$ in itself due to rotations --- i.e., translations in the second coordinate of 
$\widehat{\Omega}_{\rm pol}$. To address this issue, we extend $s_{\mu}$ to $(0,1) \times \mathbb{R}$ as follows:
$$
s_{\mu}^{\rm ext}({x})
\, =\,
s_{\mu}(x_1,  \texttt{mod}(x_2+0.5,1) - 0.5  ).
$$
}

\section{Numerical results}
\label{sec:numerics}
\subsection{Linear heat transfer in an annulus}
\label{sec:heat_transfer} 

\begin{figure}[h!]
\centering
 \subfloat[ ] 
{  \includegraphics[width=0.45\textwidth]
 {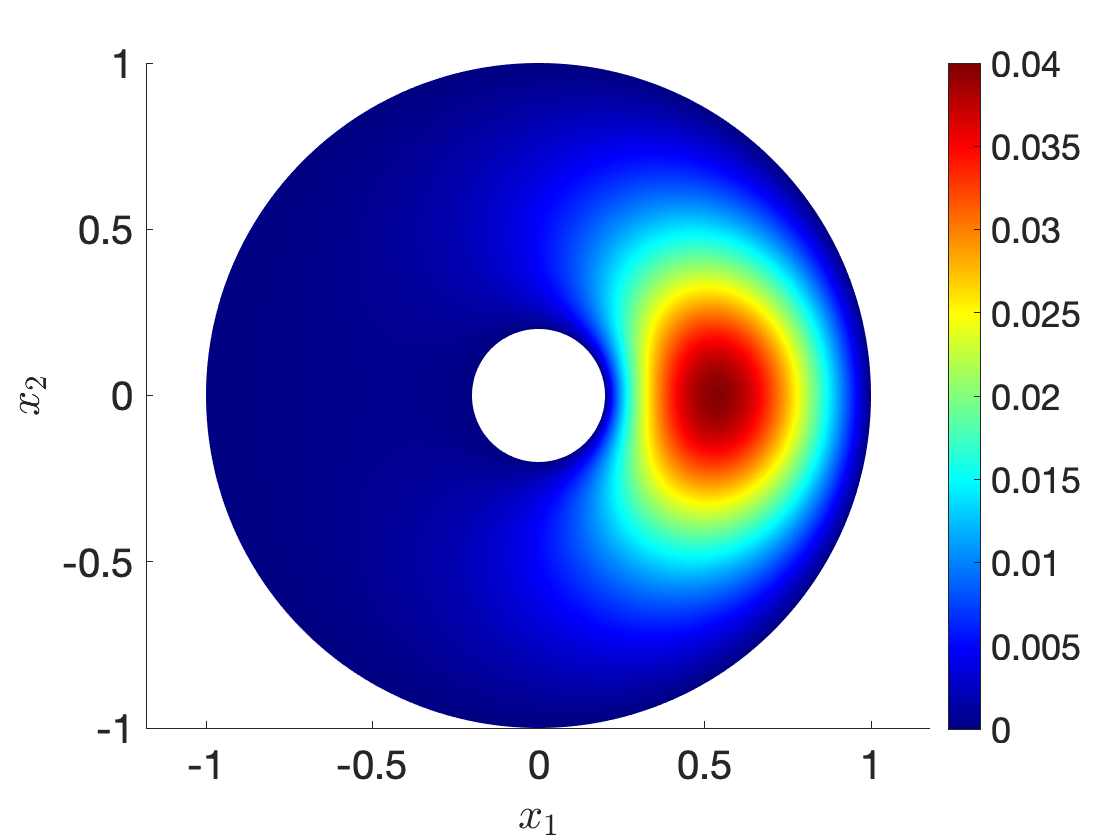}}
   ~~
 \subfloat[ ] 
{  \includegraphics[width=0.45\textwidth]
 {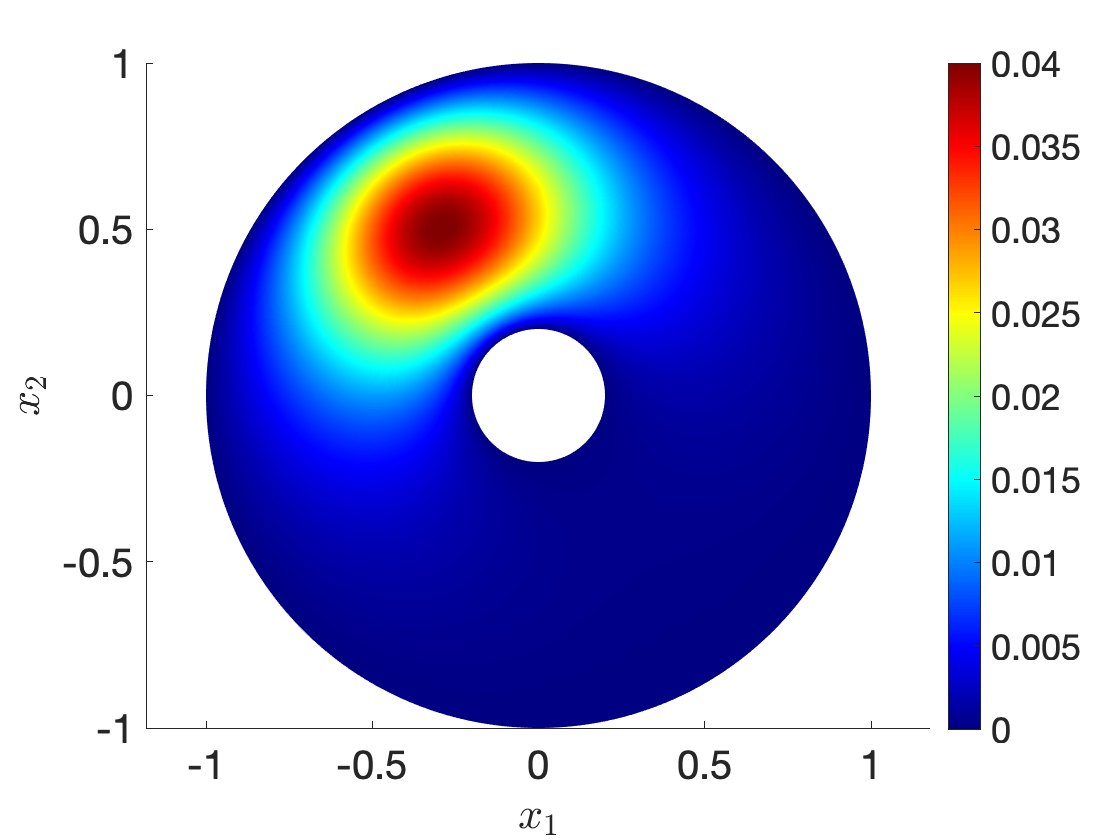}}
  
  \subfloat[ ] 
{  \includegraphics[width=0.45\textwidth]
 {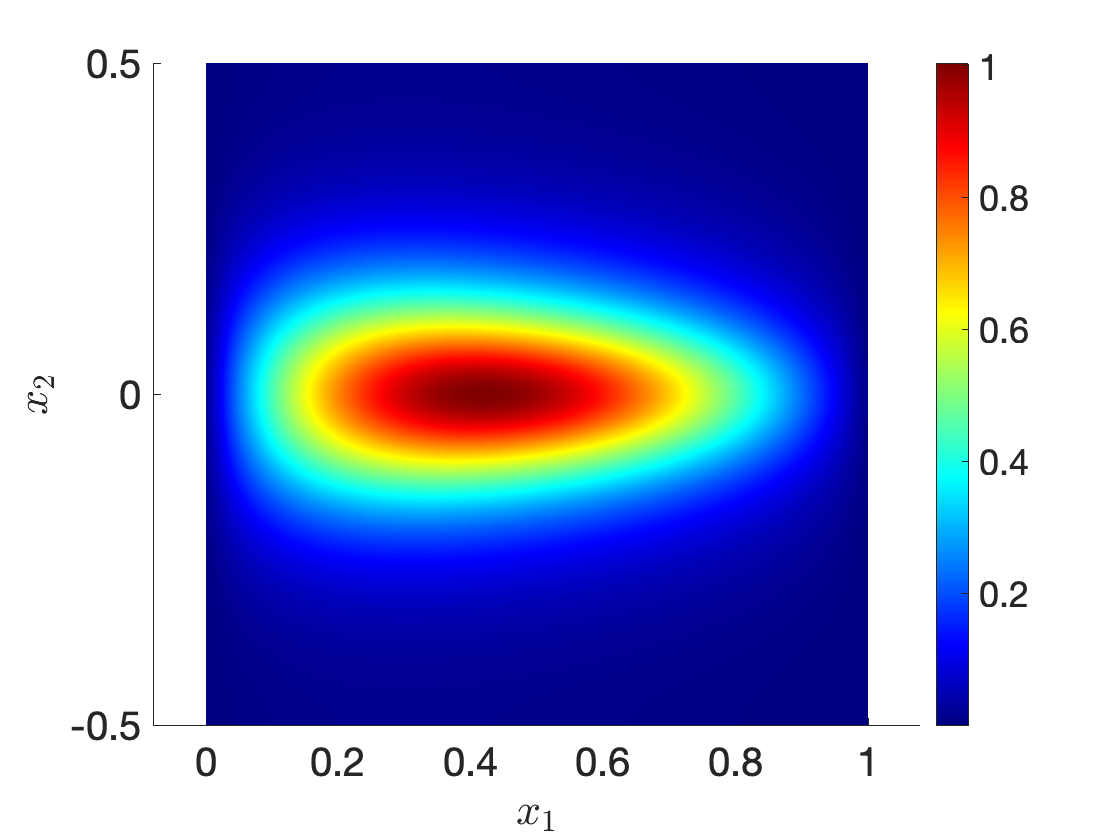}}
   ~~
 \subfloat[ ] 
{  \includegraphics[width=0.45\textwidth]
 {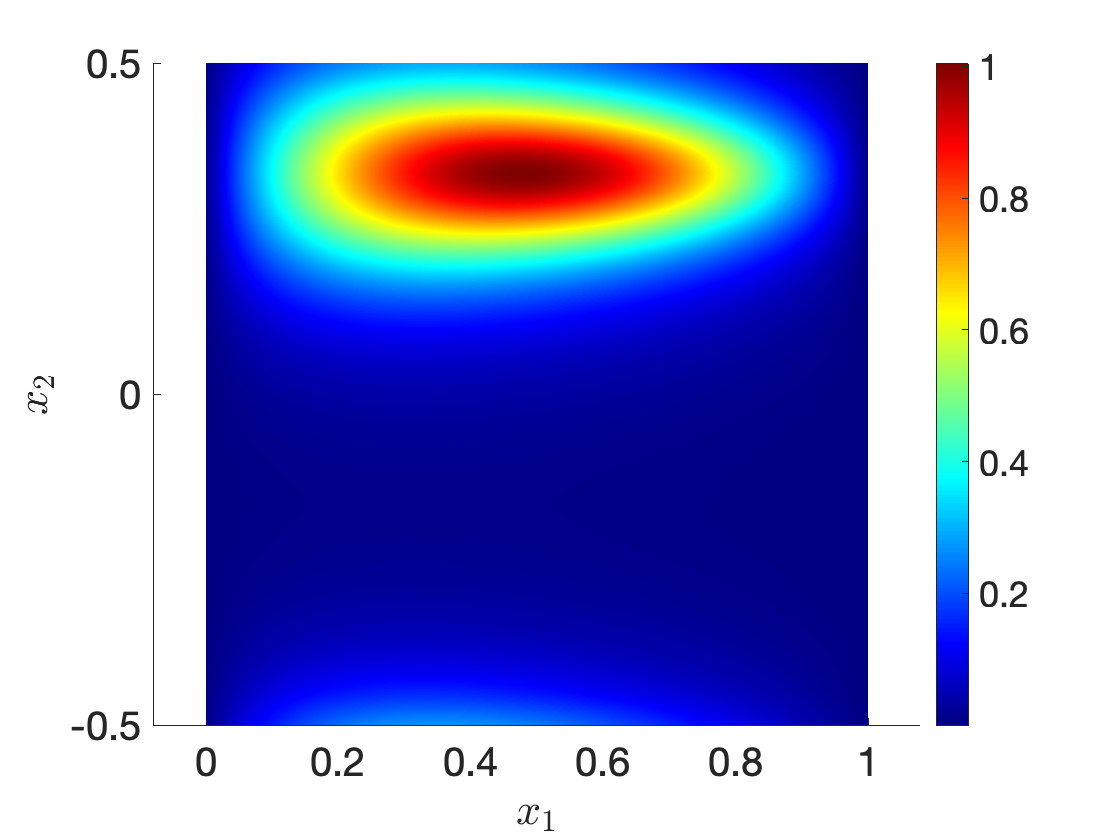}} 
  
\caption{Linear heat transfer  in an annulus.
(a) -(b) solution field in $\Omega$ for $\mu = [0,0]$ and 
 $\mu = [0.3326, 1]$.
 (c)-(d)  sensor field for the same two values of the parameter.
}
 \label{fig:annulus_vis}
  \end{figure}

We consider the problem
\begin{equation}
\label{eq:annulus_test}
\left\{
\begin{array}{ll}
- \nabla \cdot ( \kappa \nabla  u_{\mu}  ) \, = \, f_{\mu} &
{\rm in} \; \Omega = \mathcal{B}_{R=1}({0}) \setminus
\mathcal{B}_{r=0.2}({0}), \\[3mm]
u_{\mu} = 0  &
{\rm on} \; \partial \Omega;   \\
\end{array}
\right.
\end{equation}
where 
$
f_{\mu}({x}) = {\rm exp} 
\left( - 10  \| {x} - {x}_{\mu}^{\rm c}  \|_2^2   \right)
$, with 
${x}_{\mu}^{\rm c} =
 (0.5+0.1 \mu_2) 
 [\cos(2 \pi \mu_1), \sin(2 \pi \mu_1) ]$, 
$\mu = [\mu_1,\mu_2] \in \mathcal{P} =[0,1) \times [0,1]$,
$\kappa ({x}) = 0.01 + {\rm exp} \left(  - 10 | |x_1|  - 0.2 |  \right)$.
We resort to a P3 FE mesh $\mathcal{T}_{\rm hf}$ with $N_{\rm hf}=18513$ degrees of freedom. 
Figures \ref{fig:annulus_vis}(a) and (b)
show the solution field for two values of the parameter.

{We consider the sensor $s_{\mu}$ obtained using \eqref{eq:sensor_approach1} with
$\xi_{\rm s}  = 10^{-5}$ and $\widehat{N}_{\rm hf} = 3364$.
We rescale $s$ so that 
$\inf_{\widehat{\Omega}}  s_{\mu} = 0$ and 
$\sup_{\widehat{\Omega}}  s_{\mu} = 1$ for all $\mu \in  \mathcal{P}$: this simplifies the choice of the parameters $\xi, \xi_{\rm msh}$.
Figures \ref{fig:annulus_vis}(c) and (d) 
show the behavior of the sensor $s_{\mu}$ for two values of the parameter.
}

For this specific test case, we might also consider  directly a structured mesh over ${\Omega}$,   define its pre-image over $\widehat{\Omega}_{\rm pol}$ as mesh for the sensor,  and then set
$s_{\mu} = u_{\mu}^{\rm hf} \circ {\Psi}$. For more challenging problems, the use of unstructured meshes in the physical domain might be necessary: for this reason, we here choose to assess the performance of our method based on two independent meshes over $\Omega$ and $\widehat{\Omega}_{\rm pol}$.

We apply the registration procedure based on $n_{\rm train} = 10^2$ equispaced  training parameters; we set $\xi=10^{-4}$ and $\xi_{\rm msh}=10^{-6}$ in the registration statement; furthermore, we set $\bar{\mu}=[1/2,1/2]$, 
$\mathcal{S}_{N_0=1} = {\rm span} \{ s_{\bar{\mu}} \}$,
and $tol_{\rm pod}=10^{-3}$, $N_{\rm max}=5$ in Algorithm \ref{alg:registration}. The ambient space for the mapping $\mathcal{W}_{\rm hf}^{\rm pol}$ is based on \eqref{eq:ambient_space_polar} with $J_{\rm r}=12$ and $J_{\rm f} = 8$. The registration algorithm returns an expansion with  $M=2$ modes. 

We resort to RBF approximation to build the parametric mapping (cf. \eqref{eq:parametric_mapping_Phi}): as in \cite[section 4.1]{taddei2021space},  we apply RBF to each coefficient separately; we assess the goodness-of-fit using the R-squared indicator and we retain the coefficients with   R-squared larger than $0.75$. This leads to an expansion with $M=2$ modes.

\begin{figure}[h!]
\centering
 \subfloat[ ] {
\begin{tikzpicture}[scale=0.7]
\begin{loglogaxis}[
xlabel = {\LARGE {$N$}},
  ylabel = {\LARGE {$\lambda_N\lambda_1$}},
 legend entries = {unregistered,registered},
  line width=1.2pt,
  mark size=3.0pt,
  ymin=0.0000000000001,   ymax=1,
legend style={at={(0.45,0.85)},anchor=west,font=\Large}
  ]
  \addplot[line width=1.pt,color=violet,mark=square]  table {data/annulus/POD_linear.dat};
  \addplot[line width=1.pt,color=red,mark=diamond]  table {data/annulus/POD_nonlinear.dat}; 
 \end{loglogaxis}
\end{tikzpicture}
}
   ~~
 \subfloat[ ] 
{
\begin{tikzpicture}[scale=0.7]
\begin{semilogyaxis}[
xlabel = {\LARGE {$N$}},
  ylabel = {\LARGE {$E_{\rm avg}$}},
 legend entries = {unregistered,registered},
  line width=1.2pt,
  mark size=3.0pt,
  ymin=0.00001,   ymax=2,
legend style={at={(0.45,0.12)},anchor=west,font=\Large}
  ]
  \addplot[line width=1.pt,color=violet,mark=square]  table {data/annulus/error_linear.dat};
  \addplot[line width=1.pt,color=red,mark=diamond]  table {data/annulus/error_nonlinear.dat}; 
 \end{semilogyaxis}
\end{tikzpicture}
}

  \subfloat[ ] 
{  \includegraphics[width=0.45\textwidth]
 {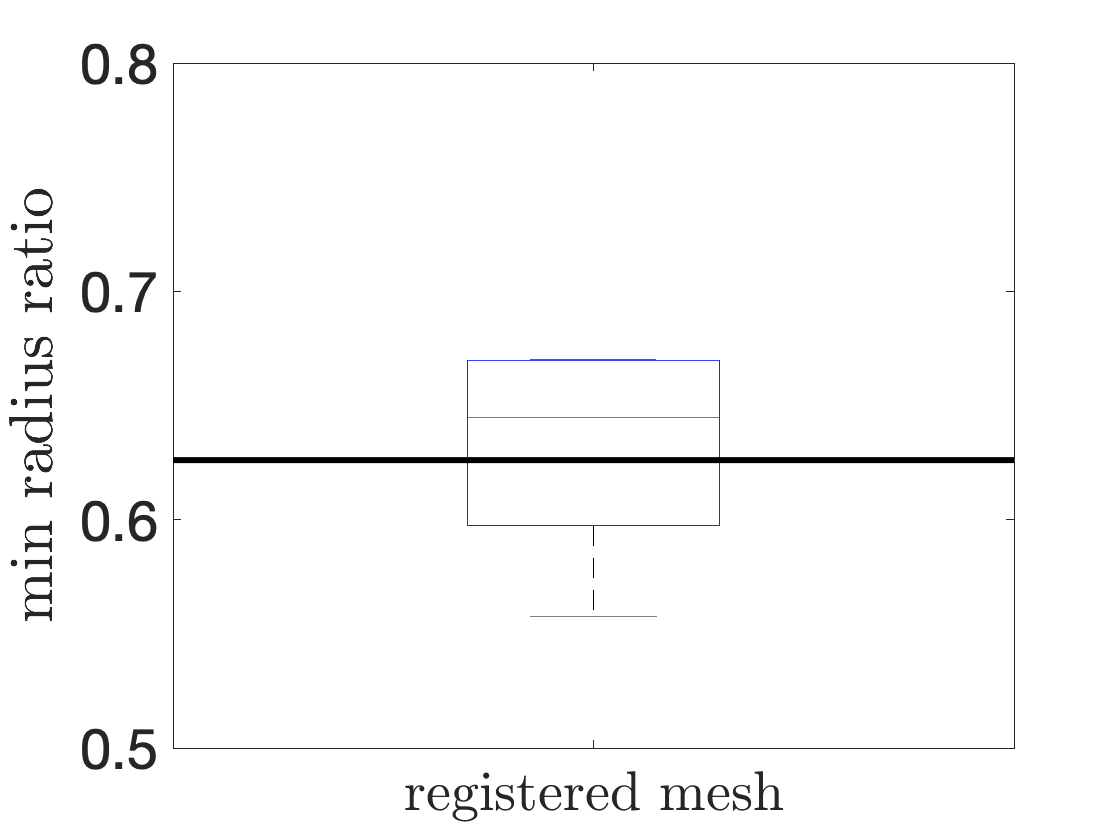}}

\caption{Linear heat transfer  in an annulus; performance.
(a)  behavior of the POD eigenvalues.
(b)  relative $H^1$ error $E_{\rm avg}$.
(c) boxplot of the minimum radius  ratio over $\mathcal{T}_{\rm hf,\mu}$.}
 \label{fig:annulus_performance}
  \end{figure}
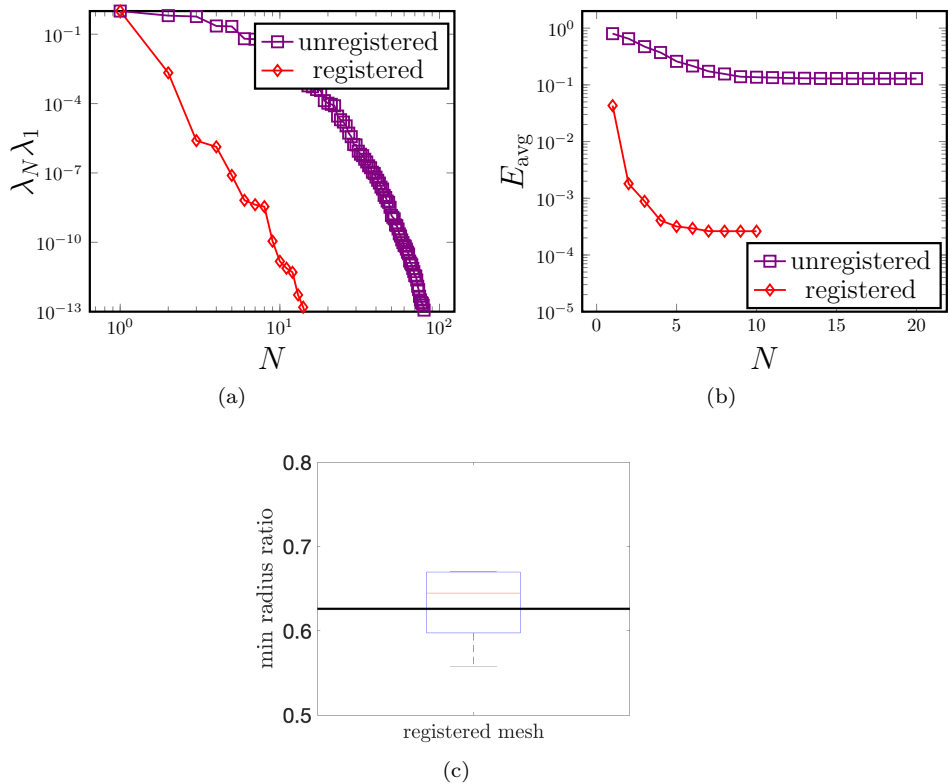

In Figure \ref{fig:annulus_performance}, we investigate online performance of our method.
In order to validate performance, we consider a set $\mathcal{P}_{\rm test} \subset \mathcal{P}$ of  $n_{\rm test}=10^2$ randomly-chosen out-of-sample parameters; {we denote by 
$E_{\rm avg}$ the average relative $H^1$ error in $\Omega$ over  the test set.}
In Figure \ref{fig:annulus_performance}(a), we show the 
behavior of the POD eigenvalues associated with the snapshots $\{  \mathbf{u}_{\mu}^{\rm hf} \}_{\mu \in \mathcal{P}_{\rm test} }$ computed using the fixed  mesh
$\mathcal{T}_{\rm hf}$ (``unregistered'') and using the parameter-dependent mesh $\mathcal{T}_{\rm hf, \mu}  = {\Phi}_{\mu}(\mathcal{T}_{\rm hf})$ (``registered'').
In Figure \ref{fig:annulus_performance}(b), we show the relative $H^1$ error $E_{\rm avg}$ for various choices of $N$ in the registered and unregistered case --- here, the reduced space and the solution coefficients  are computed using POD and RBF based on the training data, as described in section \ref{sec:parametric_reg}. 
We observe that registration significantly improves the decay of the POD eigenvalues and is also beneficial in terms of prediction. In Figure \ref{fig:annulus_performance}(c), we show the behavior of the minimum radius ratio over all elements of the meshes $\{ \mathcal{T}_{\rm hf, \mu}: \mu \in \mathcal{P}_{\rm test} \}$; the black horizontal line indicates the minimum radius ratio over all elements of the reference mesh $\mathcal{T}_{\rm hf}$.
We observe that our mapping does not significantly deteriorate the regularity of the FE mesh --- {in particular, our mapping does not lead to inverted elements for any value of the parameter in the test set.
}

\subsection{Potential flow past a parameterized airfoil}
\label{sec:airfoil} 
\begin{figure}[h!]
\centering
 \subfloat[ ] 
{  \includegraphics[width=0.55\textwidth]
 {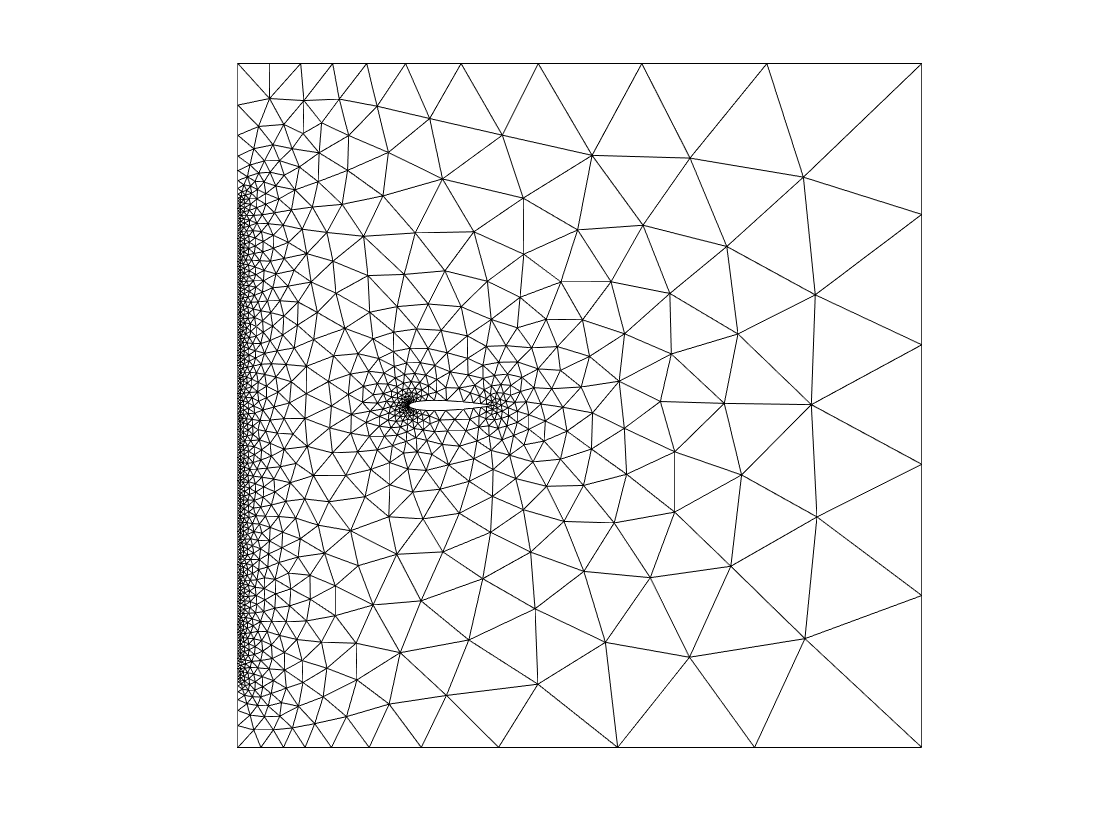}}

 \subfloat[ ] 
{  \includegraphics[width=0.45\textwidth]
 {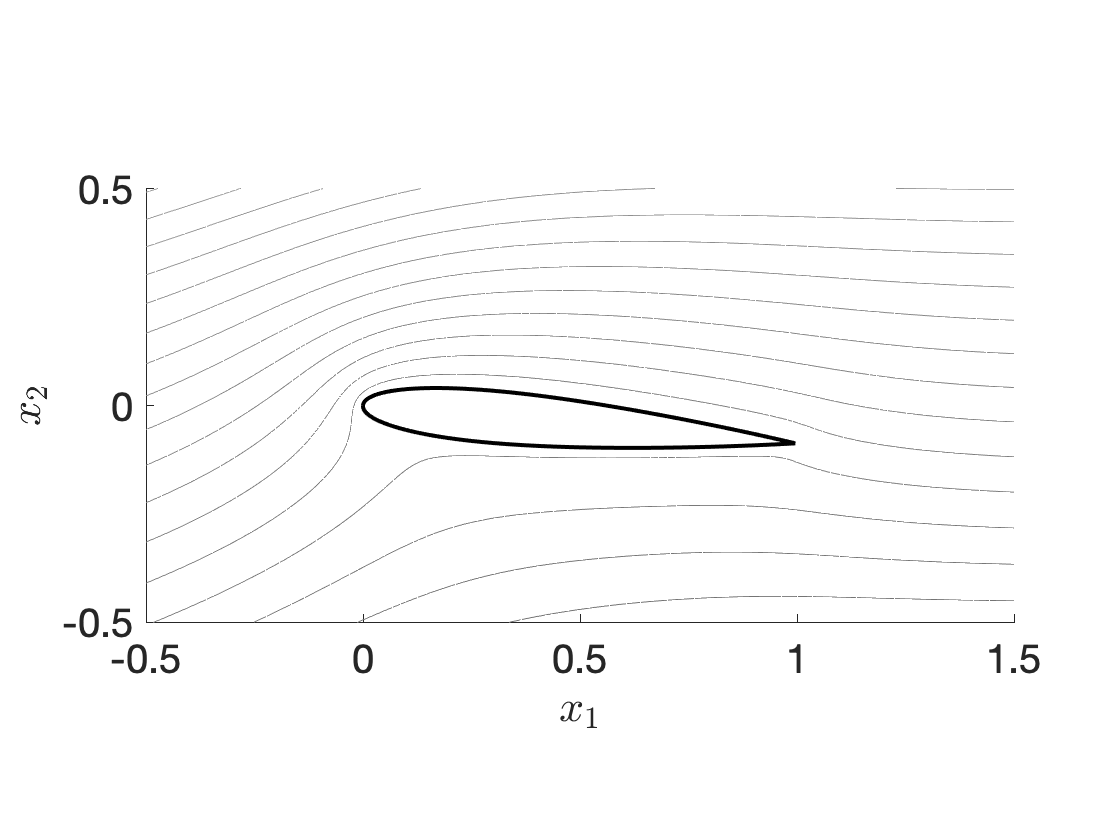}}
  ~~
  \subfloat[ ] 
{  \includegraphics[width=0.45\textwidth]
 {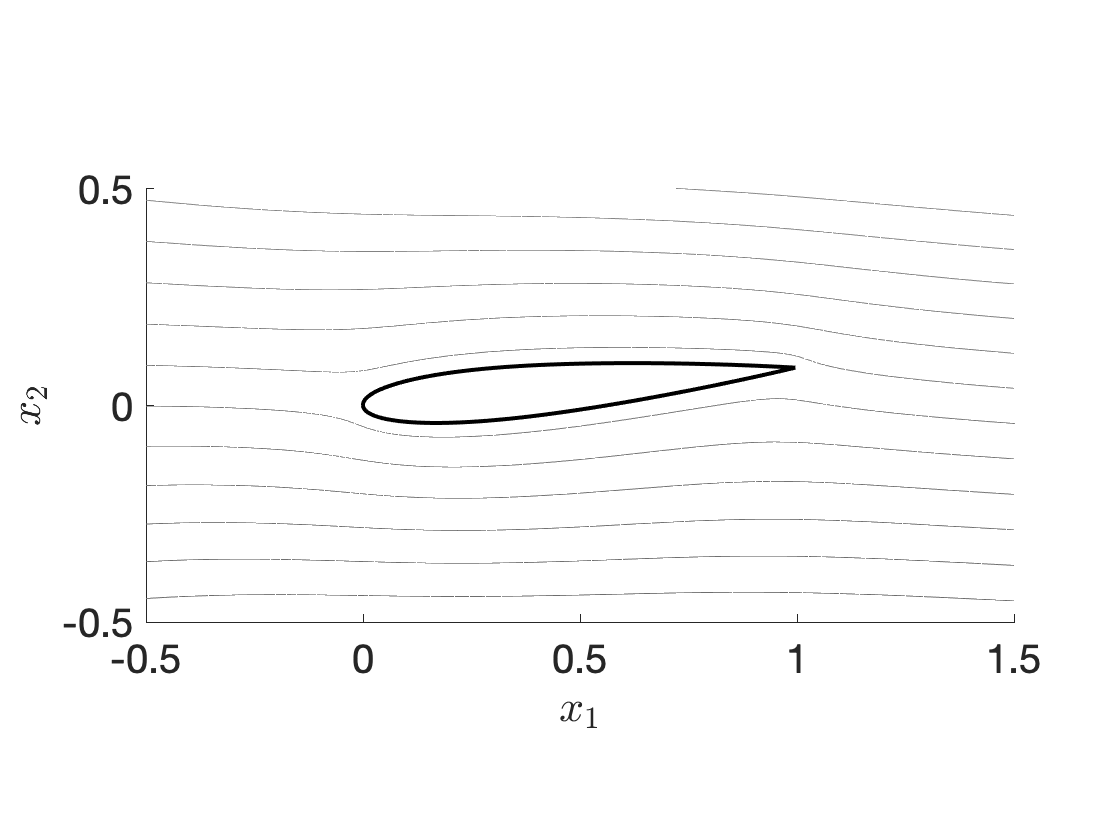}}

\caption{Potential flow past a parameterized airfoil.
(a)  computational mesh.
(b)-(c)  contour lines for $\mu=[-5/180 \pi ,   0.1   , 0.6]$ and
$\mu=[5/180 \pi ,   0.3   , 0.8]$.}
 \label{fig:naca_vis}
  \end{figure}

We consider a potential flow past a rotating airfoil.
We introduce the domain $\Omega_{\mu}= 
\Omega_{\rm box} \setminus \Omega_{\rm naca,\mu}
 \subset \mathbb{R}^2$ such that  $\Omega_{\rm box} = (x_{\rm min},x_{\rm max})\times (-H,H)$,
 $x_{\rm min} = -2, x_{\rm max} = 6, H=4$, and 
\begin{subequations}
\label{eq:naca_airfoil}
\begin{equation}
\Omega_{\rm naca,\mu} = \left\{
{\rm Rot} (\mu_3)  {x}  \; : \;
  {x} \in \Omega_{\rm naca} 
\right\},
\end{equation}
where
\begin{equation}
\begin{array}{l}
\displaystyle{
\Omega_{\rm naca} = \left\{
 {x}  \in (0,1)^2: 
 |x_2|
<f_{\rm naca,\texttt{th}}(x_1) 
\right\} 
}
\\[3mm]
\displaystyle{
f_{\rm naca,\texttt{th}}(s) = 5 \texttt{th} \left( 
0.2969 \sqrt{s} -0.1260 s - 0.3516 s^2 + 0.2843 s^3 -0.1036 s^4 \right)
},
\\
\end{array}
\end{equation}
with 
$\texttt{th}=0.12$,
and ${\rm Rot} (\theta) = [\cos(\theta), -\sin(\theta); \sin(\theta), \cos(\theta)]$ is the counterclockwise rotation matrix.
\end{subequations}
Then, we introduce the problem:
\begin{subequations}
\label{eq:potential_flow}
\begin{equation}
\left\{
\begin{array}{ll}
-\Delta u_{\mu} = 0 & {\rm in} \; \Omega_{\mu}, \\[3mm]
u_{\mu} = h_{\mu} & {\rm on} \; \partial   \Omega_{\mu}, \\
\end{array}
\right.
\end{equation}
where 
\begin{equation}
h_{\mu}({x})
\; = \;
\left\{
\begin{array}{ll}
0 &  {\rm on} \; \Gamma_{\rm btm} = (x_{\rm min}, x_{\rm max}) \times\{ -H \}, \\[3mm]
1  &  {\rm on} \; \Gamma_{\rm top} = (x_{\rm min}, x_{\rm max}) \times\{ H \}, \\[3mm]
\displaystyle{ \frac{x_2+H}{2H} }  &  {\rm on} \; \Gamma_{\rm out} = \{ x_{\rm max} \}  \times (-H, H)  , \\[3mm]
\displaystyle{ 
\frac{\bar{h}_{\mu}( \frac{x_2+H}{2H}    ) -  \bar{h}_{\mu}(0)  }{
\bar{h}_{\mu}(1) -  \bar{h}_{\mu}(0)}}  & 
 {\rm on} \; \Gamma_{\rm in} = 
\{ x_{\rm min} \}  \times (-H, H) ,       \\[3mm]
\alpha_{\rm naca, \mu}   & 
 {\rm on} \;    \partial \Omega_{\rm naca,\mu} . \\
\end{array}
\right.
\end{equation}
Here, 
$\alpha_{\rm naca, \mu} $ is chosen so that
$\nabla u_{\mu}$ is equal to zero at the trailing edge, while $ \bar{h}_{\mu}$ is given by
\begin{equation}
 \bar{h}_{\mu}(t) 
 =\frac{1}{2} \left(
1 +\frac{1}{\pi} 
\arctan \left( 100 (t - \mu_1)  \right)
 +\frac{1}{\pi} 
\arctan \left( 100 (t - \mu_2)  \right)
\right).
\end{equation}
We define the vector of parameters
$\mu=[\mu_1,\mu_2,\mu_3]$ and the parameter  domain
\begin{equation}
\mathcal{P} = [0.1,0.3] \times [0.6,0.8] \times \left[-\frac{5 \pi}{180} , \frac{5 \pi}{180}\right].
\end{equation}
We approximate the solution to  \eqref{eq:naca_airfoil} using a P3 FE discretization with  $N_{\rm hf}=10053$ degrees of freedom.
Figure \ref{fig:naca_vis} shows the computational mesh $\mathcal{T}_{\rm hf}$ in $\Omega$ and the contour lines of the solution to \eqref{eq:potential_flow} for two values of the parameter $\mu$.
\end{subequations}

We consider the partition $\{ \Omega_{q,\mu}  \}_{q=1}^{N_{\rm dd}}$, $N_{\rm dd}=4$, with  $\bar{\mu}=[0.2,0.7,0]$, depicted in Figure \ref{fig:vis_partitioned}. Recalling the definitions in section   \ref{sec:registration_sensor}, we introduce the sensor ${s}_{\mu}$
defined in \eqref{eq:sensor_approach1} based on the FE field $( \mathcal{T}_{\rm hf, \mu}^{\rm geo}, \mathbf{u}_{\mu}^{\rm  hf} )$ with $\xi_{\rm s}=10^{-4}$. We consider a 
 P3 structured grid in $\widehat{\Omega}$ with $14641$ degrees of freedom.

We apply the registration procedure based on $n_{\rm train} = 50$ randomly-sampled  training parameters; we set $\xi=10^{-4}$ and $\xi_{\rm msh}=10^{-6}$ in the registration statement; furthermore, we set $\bar{\mu}=[0.2,0.7,0]$, $\mathcal{S}_{N_0=1} = {\rm span} \{ {s}_{\bar{\mu}} \}$, $tol_{\rm pod}=10^{-3}$, $N_{\rm max}=5$ in Algorithm \ref{alg:registration}. The ambient space for the mapping $\mathcal{W}_{\rm hf,0}^{\rm dd}$ is based on \eqref{eq:calWhf0}
with $J=10$ ($M_{\rm hf} = {\rm dim} (\mathcal{W}_{\rm hf,0}^{\rm dd}) =  608$). The registration algorithm returns an expansion with  $M=10$ modes. 
Then, as in the previous example, we resort to RBF approximation to build  the parametric mapping   \eqref{eq:parametric_mapping_Phi}:
we retain an expansion with $M=5$ terms.

\begin{figure}[h!]
\centering
 \subfloat[ ] {
\begin{tikzpicture}[scale=0.7]
\begin{loglogaxis}[
xlabel = {\LARGE {$N$}},
  ylabel = {\LARGE {$\lambda_N / \lambda_1$}},
 legend entries = {unregistered,registered},
  line width=1.2pt,
  mark size=3.0pt,
  ymin=0.0000000000001,   ymax=1,
legend style={at={(0.45,0.85)},anchor=west,font=\Large}
  ]
  \addplot[line width=1.pt,color=violet,mark=square]  table {data/pflow/POD_linear.dat};
  \addplot[line width=1.pt,color=red,mark=diamond]  table {data/pflow/POD_nonlinear.dat}; 
 \end{loglogaxis}
\end{tikzpicture}
}
   ~~
 \subfloat[ ] 
{
\begin{tikzpicture}[scale=0.7]
\begin{semilogyaxis}[
xlabel = {\LARGE {$N$}},
  ylabel = {\LARGE {$E_{\rm avg}$}},
 legend entries = {unregistered,registered},
  line width=1.2pt,
  mark size=3.0pt,
  ymin=0.00001,   ymax=2,
legend style={at={(0.45,0.12)},anchor=west,font=\Large}
  ]
  \addplot[line width=1.pt,color=violet,mark=square]  table {data/pflow/error_linear.dat};
  \addplot[line width=1.pt,color=red,mark=diamond]  table {data/pflow/error_nonlinear.dat}; 
 \end{semilogyaxis}
\end{tikzpicture}
}

  \subfloat[ ] 
{  \includegraphics[width=0.5\textwidth]
 {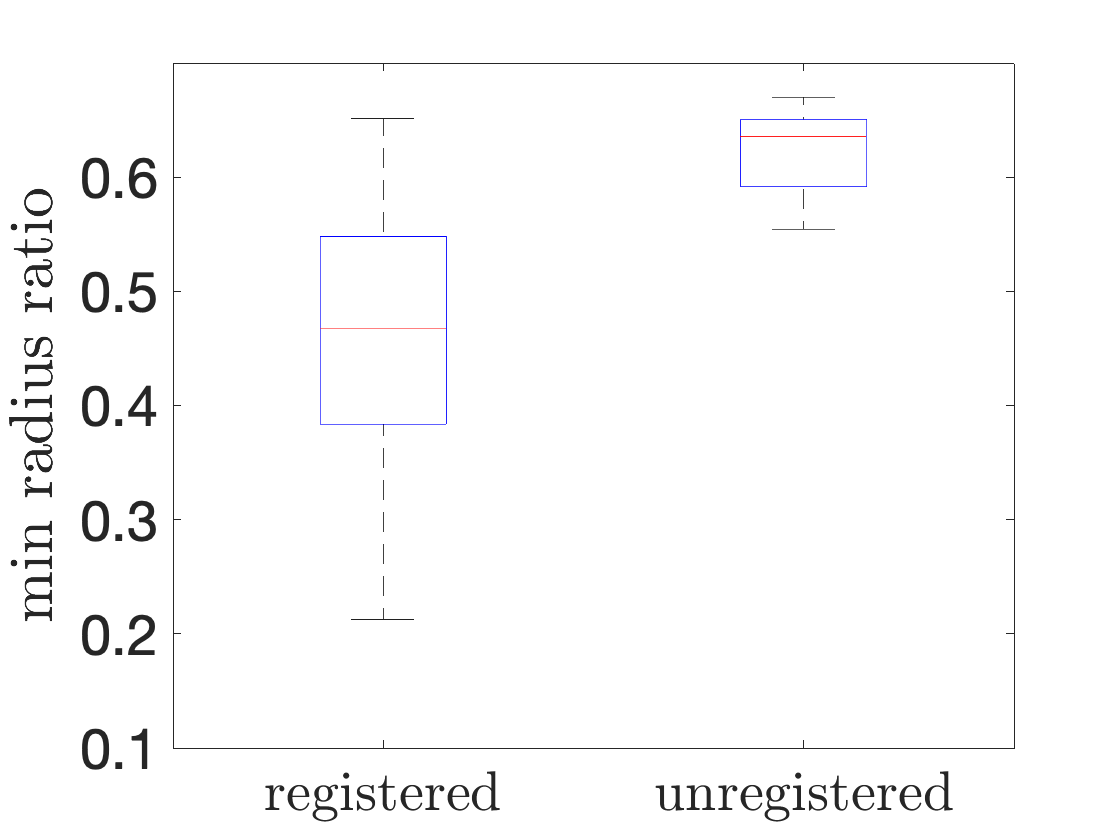}}
~~
  \subfloat[ ] 
{  \includegraphics[width=0.5\textwidth]
 {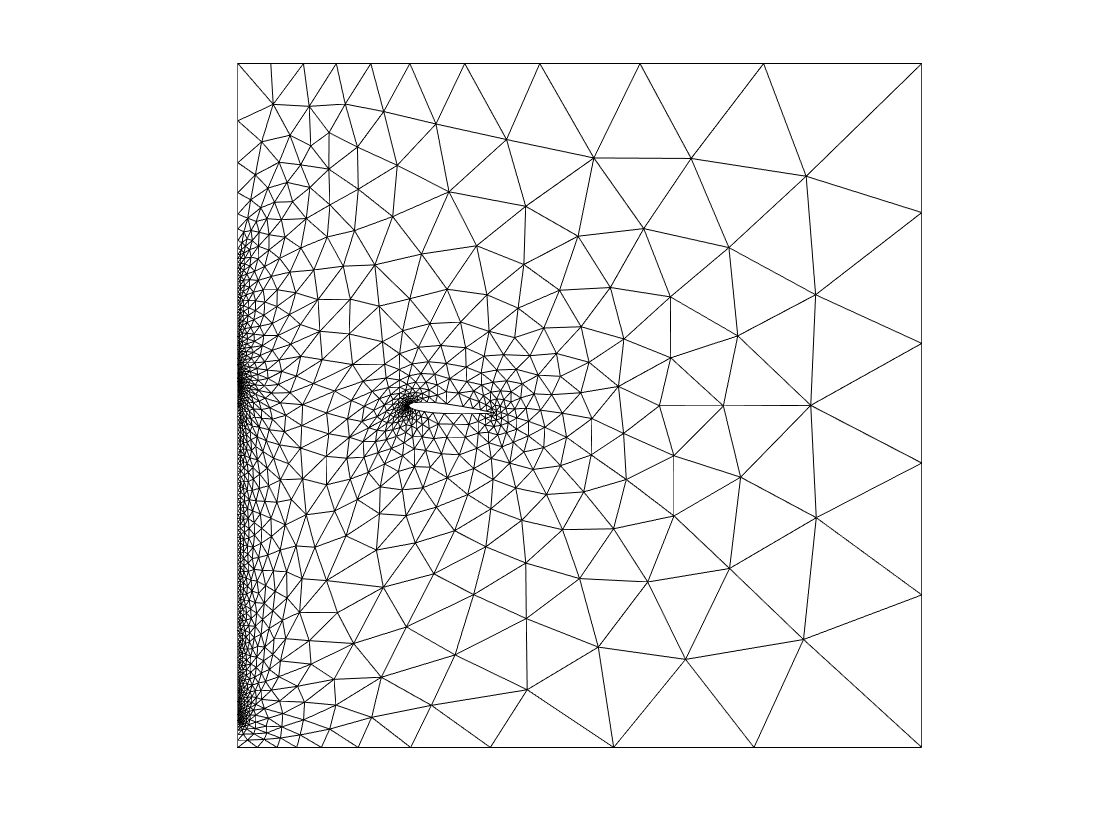}}

\caption{Potential flow past a parameterized airfoil; performance.
(a)  behavior of the POD eigenvalues.
(b)  relative $H^1$ error $E_{\rm avg}$.
(c) boxplots of the minimum radius  ratio over $\mathcal{T}_{\rm hf,\mu}$ and $\mathcal{T}_{\rm hf,\mu}^{\rm geo}$.
(d) deformed mesh $\mathcal{T}_{\rm hf,\mu}$ for $\mu=[-5/180 \pi ,   0.1   , 0.6]$.}
 \label{fig:naca_performance}
  \end{figure}
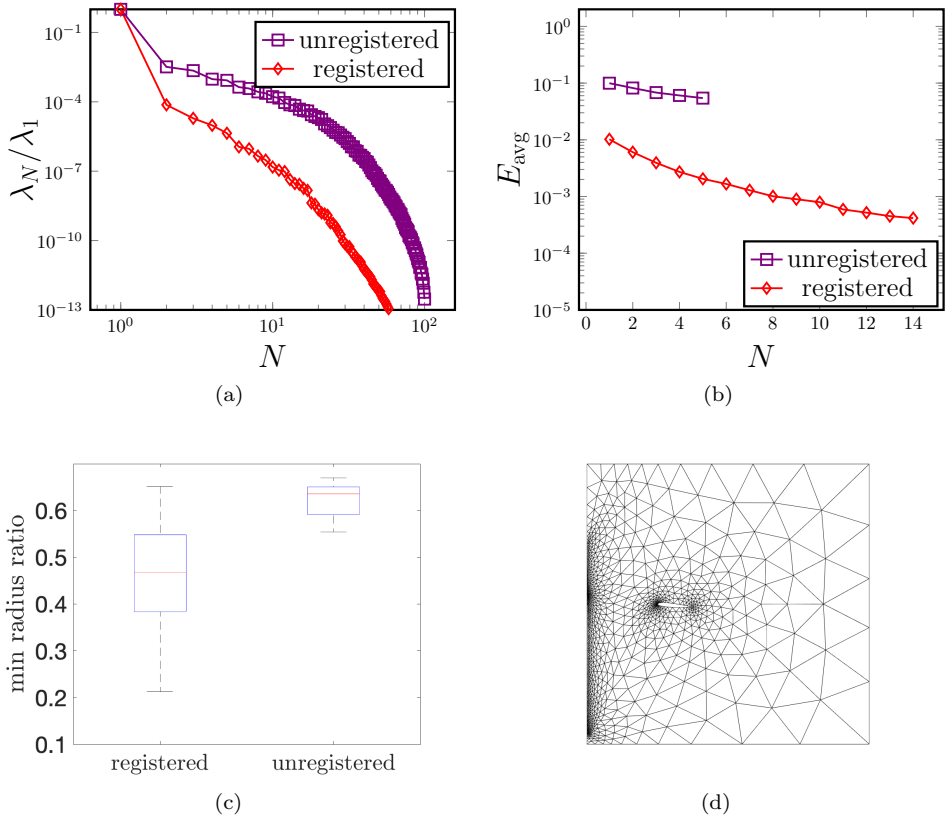

In Figure \ref{fig:naca_performance}, we investigate online performance of our method.
In order to validate performance, we consider a set $\mathcal{P}_{\rm test} \subset \mathcal{P}$ of  $n_{\rm test}=10^2$ randomly-chosen out-of-sample snapshots. 
In Figure \ref{fig:naca_performance}(a), we show the 
behavior of the POD eigenvalues associated with the snapshots $\{  \mathbf{u}_{\mu}^{\rm hf} \}_{\mu \in \mathcal{P}_{\rm test} }$ computed using the a priori deformed  mesh
$\mathcal{T}_{\rm hf,\mu}^{\rm geo}$ (``unregistered'') and using the parameter-dependent solution-aware  mesh $\mathcal{T}_{\rm hf,\mu}  = {\Phi}_{\mu}(\mathcal{T}_{\rm hf})$ (``registered'').
In Figure \ref{fig:naca_performance}(b), we show the average relative $H^1$ error $E_{\rm avg}$  on the test set for various choices of $N$ in the registered and unregistered case --- as for the mapping, we only retain the coefficients that ``pass'' the goodness-of-fit test based on the R-squared indicator. 
As in the previous test case,   registration significantly improves the decay of the POD eigenvalues and is also beneficial in terms of prediction. In Figure \ref{fig:naca_performance}(c), we show the behavior of the minimum radius ratio over all elements of the meshes $\{ \mathcal{T}_{\rm hf,\mu}: \mu \in \mathcal{P}_{\rm test} \}$ (``registered") and
$\{ \mathcal{T}_{\rm hf,\mu}^{\rm geo}: \mu \in \mathcal{P}_{\rm test} \}$
(``unregistered").
Although  the solution-aware  mapping deteriorates the regularity of the FE mesh for this test case, the minimum radius ratio remains larger than $0.2$ for all test cases --- and can be controlled by suitably changing $\mathfrak{f}_{\rm msh, max}$ in \eqref{eq:local_mesh_distortion}.
{In Figure \ref{fig:naca_performance}(d),  we show the deformed mesh $\mathcal{T}_{\rm hf,\mu}$ for
$\mu=[-5/180 \pi ,   0.1   , 0.6]$.
}

\subsection{Inviscid transonic flow past a parameterized airfoil}
\label{sec:euler_equation}

We consider the problem of approximating the solution to the  compressible Euler equations with varying Mach number, past a rotating airfoil with varying upper and lower thicknesses.  
We consider   the domain $\Omega_{\mu}= 
\Omega_{\rm box} \setminus \Omega_{\rm naca,\mu}
 \subset \mathbb{R}^2$ such that  $\Omega_{\rm box} = (-4,10)\times (-10,10)$, and
\begin{equation}
\label{eq:euler_airfoil}
\Omega_{\rm naca,\mu} = \left\{
{\rm Rot} (\mu_3)  {x}  \; : \;
x \in (0,1)^2, \;
- f_{\rm naca,\mu_1}(x_1) 
  {x}_2 < f_{\rm naca,\mu_1}(x_2) 
\right\},
\end{equation}
where $f_{\rm naca,\texttt{th}}$ is defined in \eqref{eq:naca_airfoil} and ${\rm Rot} (\mu_3)$ is the counterclockwise rotation matrix. Then, we introduce the vector of conserved variables $U = [\rho, \rho u, E]$ where $\rho$ is the fluid density, $u=[u_1,u_2]$ is the velocity, and $E$ is the total energy; we further define the pressure 
$p  = (\gamma-1) (E - 1/2 \rho \| u \|_2^2)$ with $\gamma=1.4$, the sound velocity
$a = \sqrt{\frac{\gamma p}{\rho}}$ and the 
Mach number ${\rm Ma} = \frac{\| u  \|_2}{a}$.

We define  the steady Euler equations 
\begin{equation}
\label{eq:euler_equation}
\nabla \cdot F(U_{\mu}) = 0  \; {\rm in} \; \Omega_{\mu}, \quad
{\rm where} \; 
F(U) = \left[
\begin{array}{l}
\rho u^T \\
\rho u u^T + p \mathbbm{1} \\
u^T (E+p) \\
\end{array}
\right].
\end{equation}
We impose  wall   boundary conditions on the airfoil, transmissive conditions on the lower and upper boundaries, 
we set a constant parametric inflow boundary condition such that the Mach number ${\rm Ma}$ is equal to $\mu_4$; finally, we impose 
constant pressure at the outflow equal to the inflow pressure.

\begin{figure}
\subfloat[$\mu_{\rm test}^1$] 
{  \includegraphics[width=0.45\textwidth]
 {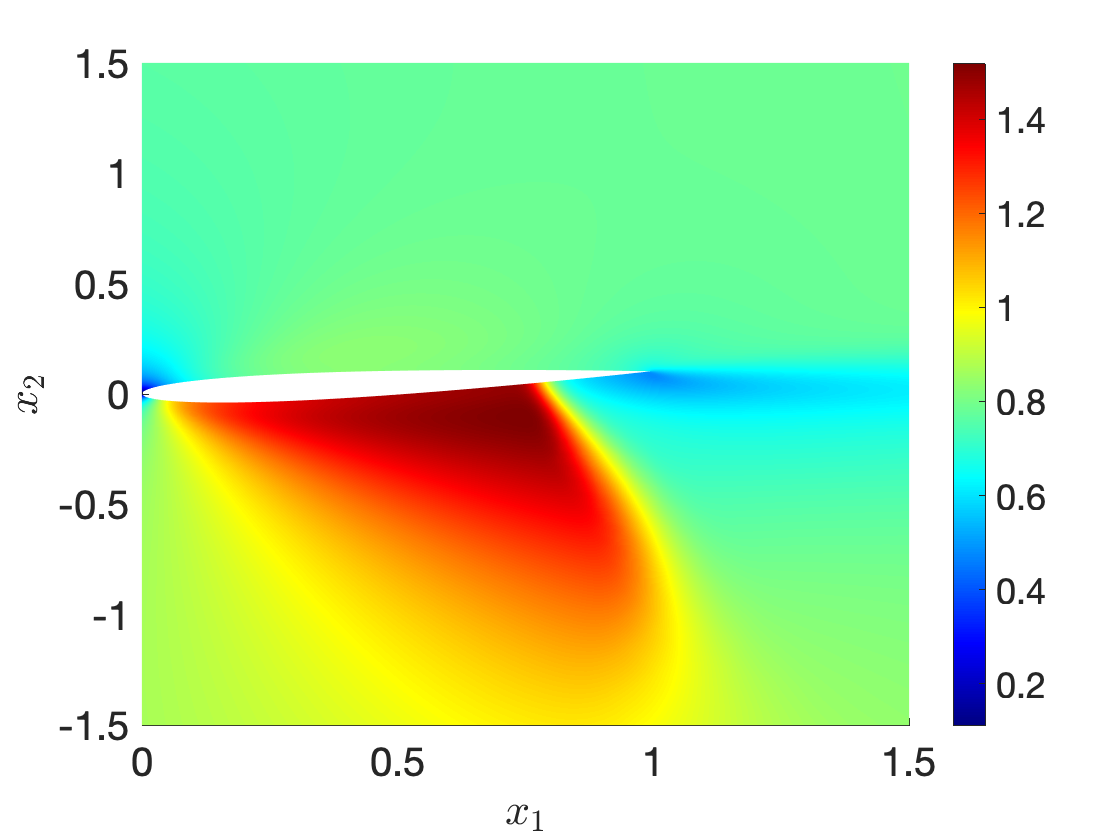}}
~~
  \subfloat[$\mu_{\rm test}^2$] 
{  \includegraphics[width=0.45\textwidth]
 {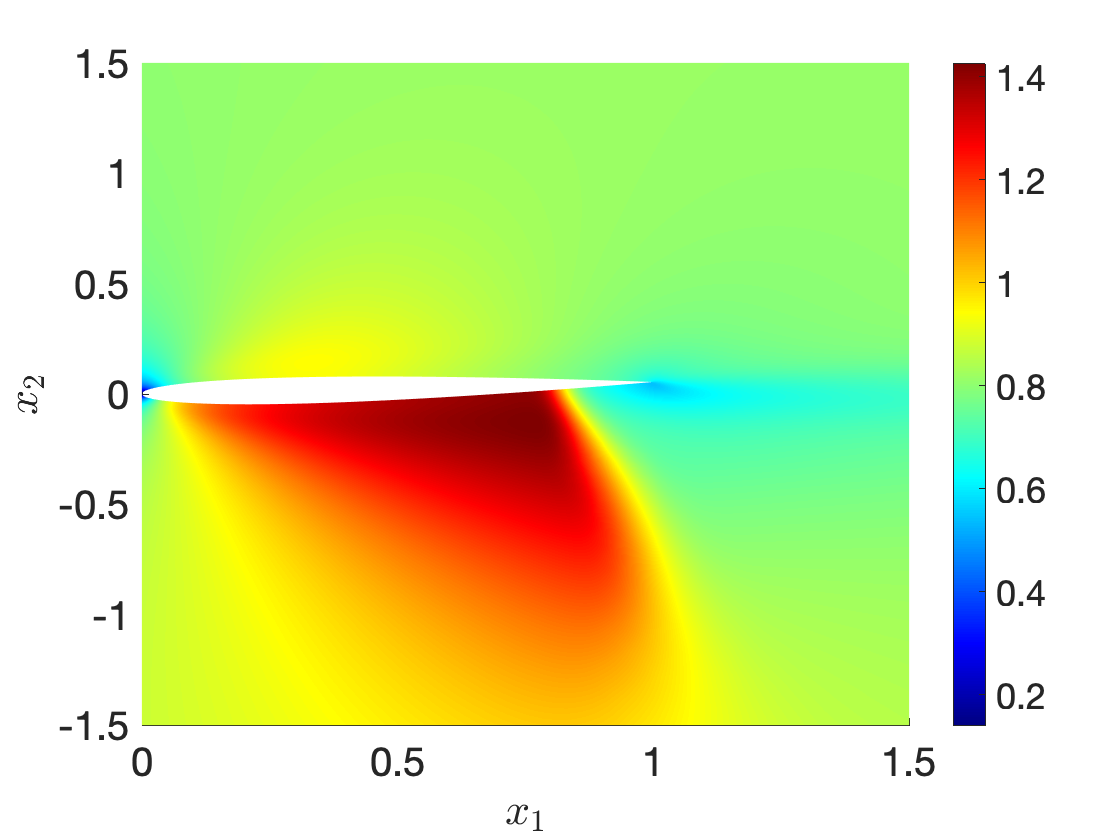}}

\caption{Transonic flow past a parameterized airfoil. 
Behavior of the Mach number for  $\mu_{\rm test}^1=[0.1226,    0.1222,    5.9072^o,    0.7953]$ and 
$\mu_{\rm test}^2= [0.1181,    0.1207,    3.1169^o,    0.8201]$.}
 \label{fig:euler_vis}
  \end{figure}

We define the vector of parameters $\mu=[\mu_1,\mu_2,\mu_3,\mu_4]$ and the parameter region:
$$
\mathcal{P} = 
\left[ 0.95 {\rm th}, 1.05 {\rm th} \right]
\times
\left[ 0.95 {\rm th}, 1.05 {\rm th} \right]
\times
\left[ 3^o, 6^o \right]
\times
[0.77,0.83],
$$
with ${\texttt{th}}=0.12$.
Note that for all parameter values the flow is transonic. We resort to a P3 discontinuous Galerkin (DG) discretization with $N_{\rm e}=9828$ elements ($N_{\rm hf}=393 120$). Figure \ref{fig:euler_vis} shows the behavior of the Mach number for two values of the parameter as predicted by our hf code. Geometry parameterization and the subsequent coarse-grained partition used for registration are the same used in section \ref{sec:airfoil}.

We train our ROM based on $n_{\rm train}=50$ randomly-sampled parameters and we 
assess performance based on $n_{\rm test}=10$ out-of-sample configurations. 
Registration sensor is computed based on the Mach number using the approach in \eqref{eq:sensor_approach2} with $\xi_{\rm s}=10^{-4}$. Algorithm \ref{alg:registration} is then applied using the same specifications as in section \ref{sec:airfoil}. The resulting map has $M=3$ terms. Computation of the solution coefficients associated with the POD modes is performed using RBF regression; we set to zero coefficients that do not pass the goodness-of-fit test.
We measure performance in terms of the average relative $L^2(\Omega_{\mu})$ error on the test set.

\begin{figure}[h!]
\centering
 \subfloat[ ] {
\begin{tikzpicture}[scale=0.7]
\begin{loglogaxis}[
xlabel = {\LARGE {$N$}},
  ylabel = {\LARGE {$E_{\rm avg}^{\rm proj}$}},
 legend entries = {unregistered,registered},
  line width=1.2pt,
  mark size=3.0pt,
  ymin=0.00001,   ymax=0.1,
legend style={at={(0.45,0.85)},anchor=west,font=\Large}
  ]
  \addplot[line width=1.pt,color=violet,mark=square]  table {data/euler/proj_linear.dat};
  \addplot[line width=1.pt,color=red,mark=diamond]  table {data/euler/proj_nonlinear.dat}; 
 \end{loglogaxis}
\end{tikzpicture}
}
   ~~
 \subfloat[ ] 
{
\begin{tikzpicture}[scale=0.7]
\begin{semilogyaxis}[
xlabel = {\LARGE {$N$}},
  ylabel = {\LARGE {$E_{\rm avg}$}},
 legend entries = {unregistered,registered},
  line width=1.2pt,
  mark size=3.0pt,
  ymin=0.00001,   ymax=0.1,
legend style={at={(0.45,0.85)},anchor=west,font=\Large}
  ]
  \addplot[line width=1.pt,color=violet,mark=square]  table {data/euler/error_linear.dat};
  \addplot[line width=1.pt,color=red,mark=diamond]  table {data/euler/error_nonlinear.dat}; 
 \end{semilogyaxis}
\end{tikzpicture}
}

\caption{Transonic flow past a parameterized airfoil: performance of registration.
(a)  behavior of the relative projection error.
(b)  relative $L^2$ error $E_{\rm avg}$.}
 \label{fig:euler_performance}
  \end{figure}
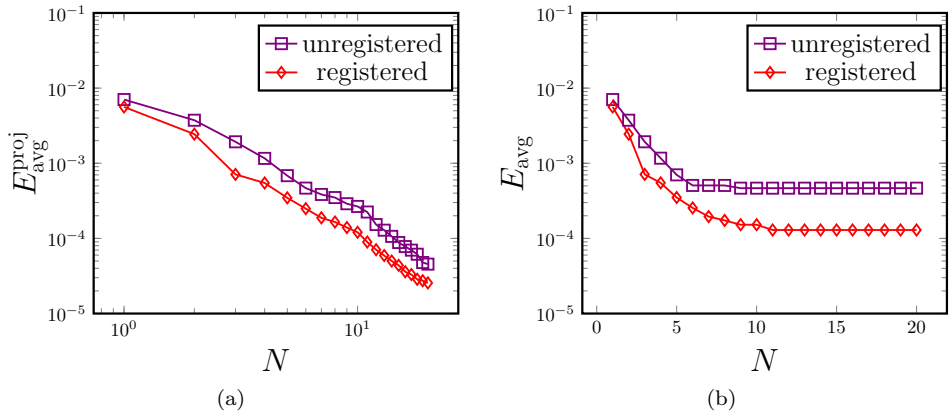

\begin{figure}
\subfloat[unregistered] 
{  \includegraphics[width=0.45\textwidth]
 {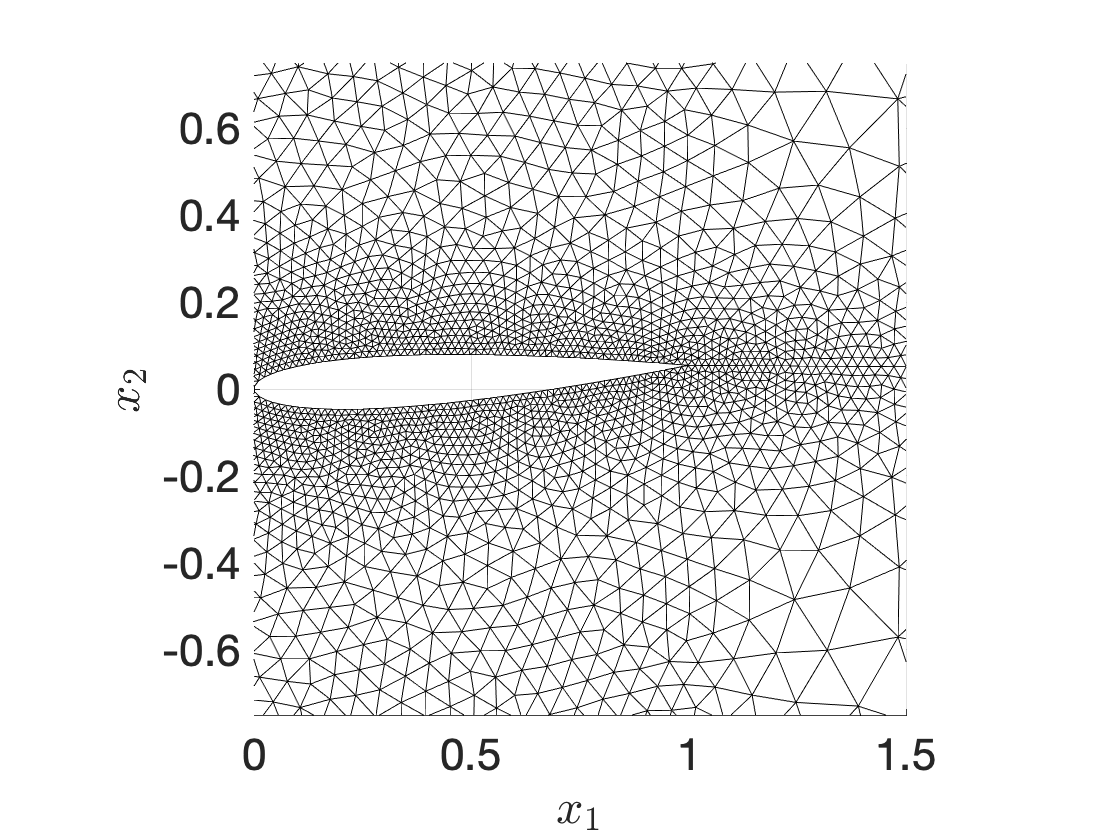}}
~~
  \subfloat[unregistered] 
{  \includegraphics[width=0.45\textwidth]
 {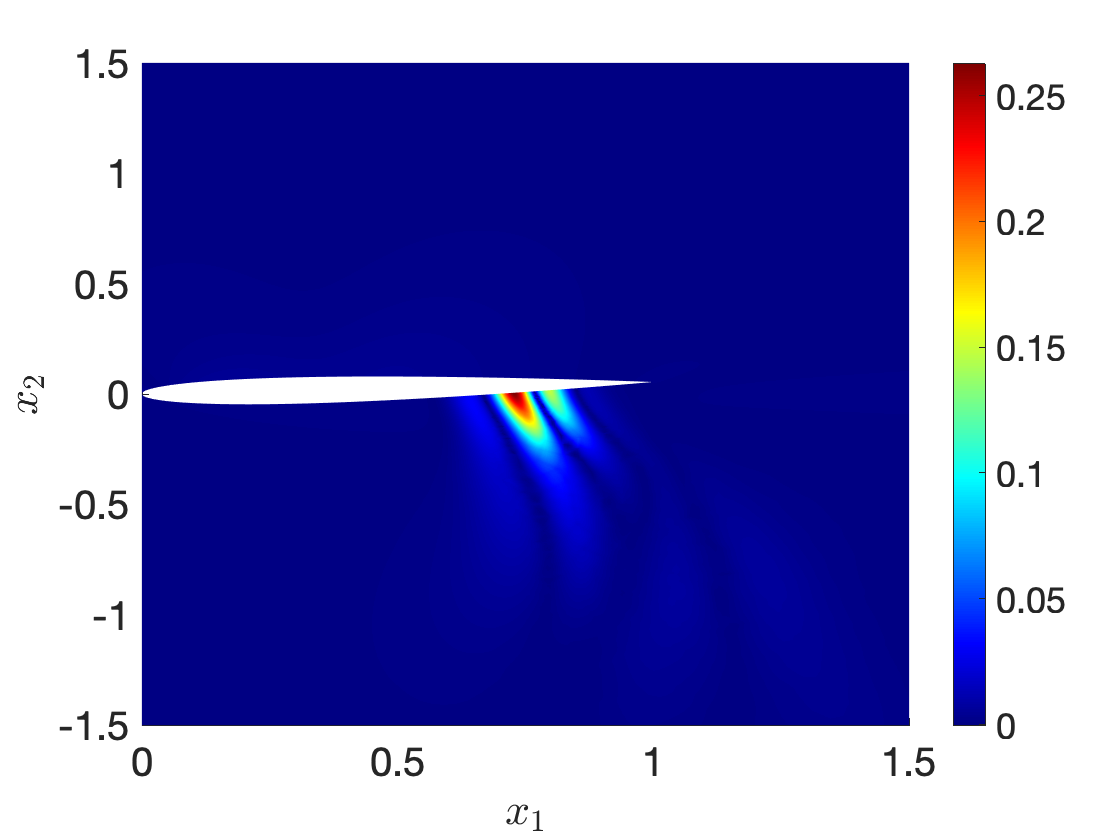}}

\subfloat[registered] 
{  \includegraphics[width=0.45\textwidth]
 {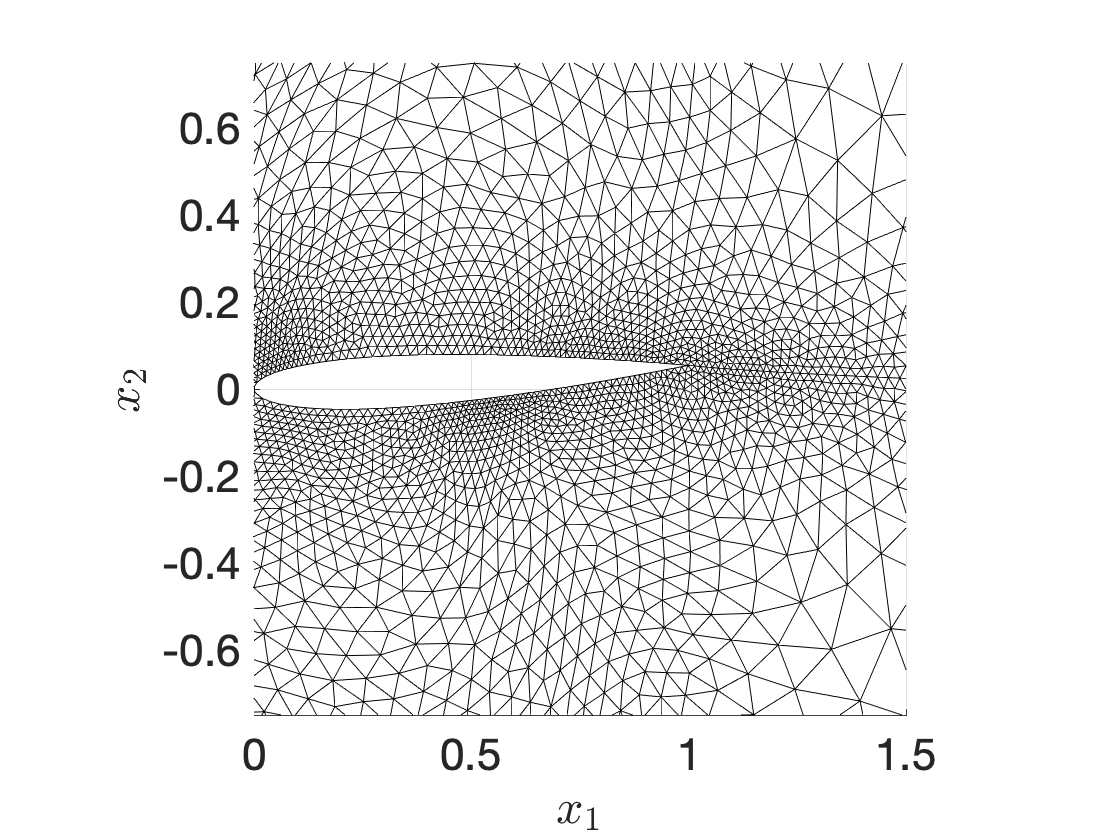}}
~~
  \subfloat[registered] 
{  \includegraphics[width=0.45\textwidth]
 {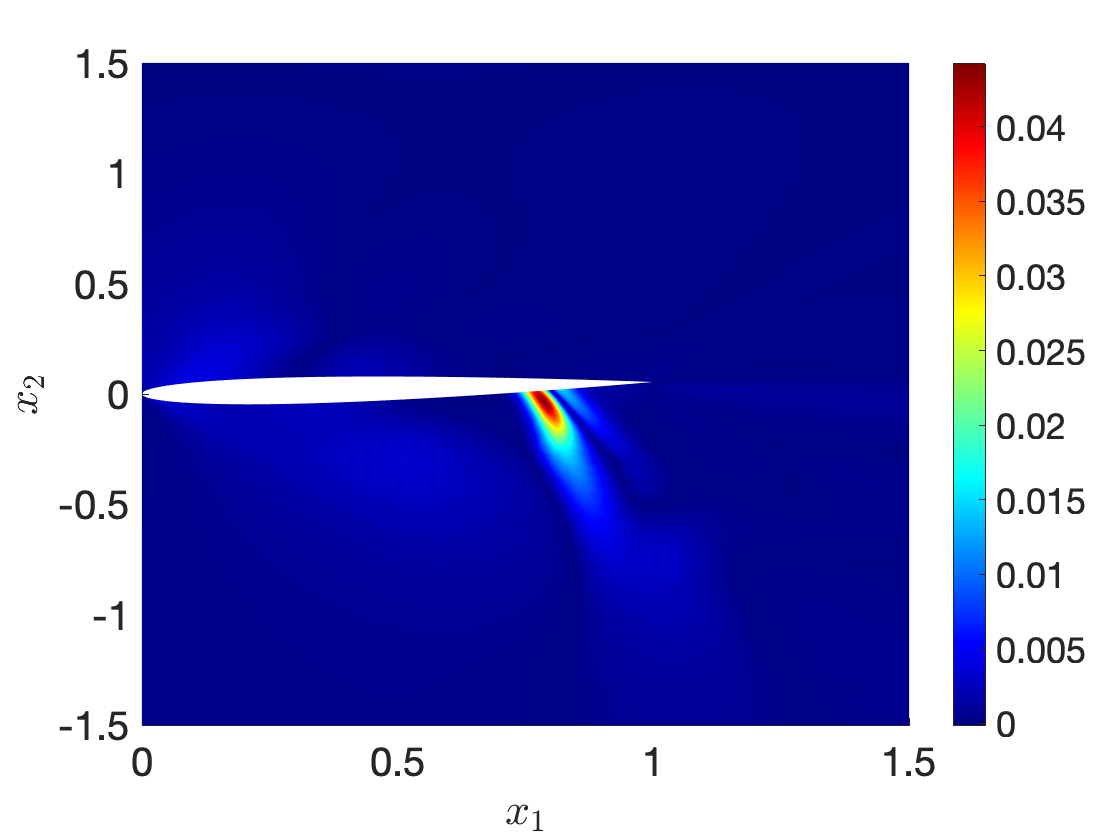}}

\caption{Transonic flow past a parameterized airfoil: performance of registration.
Comparison between unregistered and registered approaches for $\mu_{\rm test}^2=[0.1181,    0.1207,    3.1169^o,    0.8201]$.
(a)-(c) FE meshes.
(b)-(d) behavior of the prediction error in Mach estimation.}
 \label{fig:euler_vis_more}
  \end{figure}

Figure \ref{fig:euler_performance} shows the behavior of the relative projection error and of the relative prediction error based on RBF approximation for the registered and the unregistered methods. As in the previous examples, registration significantly helps improve performance for out-of-sample configurations.
Figure \ref{fig:euler_vis_more} shows the deformed meshes $\mathcal{T}_{\rm hf, \mu}^{\rm geo}$ and $\mathcal{T}_{\rm hf, \mu}= \Phi_{\mu}( \mathcal{T}_{\rm hf})$ and the error 
between hf and predicted Mach number for registered and unregistered ROMs, for 
 $\mu_{\rm test}^2=[0.1181,    0.1207,    3.1169^o,    0.8201]$. We observe that the registered ROM offers significantly more accurate performance, especially in the proximity of the shock.

\section{Summary and discussion}
\label{sec:conclusions}

In this paper, 
we extended the  registration algorithm first introduced in \cite{taddei2020registration,taddei2021space} to geometries $\Omega$ that are not isomorphic to the unit square. 
We devised a specialized approach for annular domains based on a  spectral expansion in a reference domain; then, we introduced a partitioned approach for general two-dimensional geometries based on a suitable spectral element approximation.
Results for three model problems --- a heat-transfer problem with a moving source, a potential flow past a rotating airfoil with parameterized inflow condition, and a transonic flow past a rotating non-symmetric airfoil with varying thickness --- demonstrate the effectiveness of the approach.

{However, the registration approach is still limited to relatively simple geometries with smooth boundaries. The extension to more challenging problems, which arise for instance in coastal engineering and  hydraulics,
does require major advances and 
 is the subject of ongoing research.
 An alternative route, which we are also currently investigating, is to apply the present approach away from non-smooth boundaries $\Gamma_{\rm co}$:  this can be achieved by setting $\Phi=\texttt{id}$ in a neighborhood of $\Gamma_{\rm co}$.
}

As discussed in section \ref{sec:parametric_reg}, a major limitation of our approach is the need for several offline simulations to build the mapping ${\Phi}$.
{The need for sufficiently-dense discretizations of the parameter domain might indeed preclude the application of our method to high-dimensional parameterizations.
}
 To address this issue, we wish  to devise multi-fidelity strategies to build the dataset of registration sensors $\{ {s}_{\mu^k} \}_k$ and then resort to standard greedy algorithms to devise the linear approximation $\mu \mapsto \widehat{\mathbf{u}}_{\mu}$.
{More in detail, since registration is performed before applying model reduction and the snapshots used for registration are not directly employed for the subsequent construction of the ROM (cf. section \ref{sec:parametric_reg}), we can readily  combine different sources of information of different fidelity to generate the dataset of sensors. 
}
 
Alternatively,  we wish to devise a fully-intrusive  projection-based  technique for the simultaneous computation of solution and mapping coefficients. In this respect, the work by Zahr and Persson \cite{zahr2018optimization} in the discontinuous Galerkin finite element  framework might provide the foundations for our pMOR procedure.

Furthermore, we wish to investigate what phenomena of interest in science and engineering can be properly described using registration-based methods: this question requires a deep understanding of the physical and mathematical properties of the system of interest,  and  clearly lies at the intersection between mechanics and applied mathematics.

\section*{Acknowledgements}
The authors acknowledge the support by European Union’s Horizon 2020 research and innovation programme under the Marie Skłodowska-Curie Actions, grant agreement 872442 (ARIA).
The authors thank  Professor Angelo Iollo (Inria Bordeaux), 
Dr. Andrea Ferrero (Politecnico di Torino), 
 Dr. Cédric Goeury and Dr.  Angélique Ponçot (EDF) for fruitful discussions.
 
\bibliographystyle{plain}
\bibliography{all_refs}
 
\end{document}